\documentclass{article}
\usepackage{latexsym,amsmath,amssymb,amscd, bm}
\usepackage[mathscr]{eucal}
\usepackage[dvips]{color}
\begin{document}

\title{On a Mathematical Model of the Rotating Atmosphere of the Earth  }
\author{Tetu Makino \footnote{Professor Emeritus at Yamaguchi University, Japan;  E-mail: makino@yamaguchi-u.ac.jp} }
\date{\today}
\maketitle

\newtheorem{Lemma}{Lemma}
\newtheorem{Proposition}{Proposition}
\newtheorem{Theorem}{Theorem}
\newtheorem{Definition}{Definition}
\newtheorem{Remark}{Remark}
\newtheorem{Corollary}{Corollary}
\newtheorem{Notation}{Notation}
\newtheorem{Assumption}{Assumption}
\newtheorem{Approximation}{Approximation}
\newtheorem{OpenProposition}{Open Proposition}

\numberwithin{equation}{section}

\begin{abstract}
In meteorology the analysis of motions of the atmosphere on the Earth has been done using various mathematical models and using various approximations. In this article as the simplest model the compressible Euler equations with barotropic equation of state of the ideal gas is analyzed under the co-ordinate system which rotates with constant angular velocity. Mathematically rigorous inquiry is tried. Although problems remain to be open, some fundamental results are exhibited.



\end{abstract}

\section{Introduction}

In meteorology the analysis of motions of the atmosphere on the Earth has been done using various mathematical models and using various approximations. Actually models used for numerical weather prediction or climate simulation must take into account various factors besides the fundamental state variables of the gas and must meet practical computational efficiency with moderate accuracy
in the process of the numerical simulations. Discussions of choice of models have been accumulated as an enormous collection in the meteorological literatures. However fundamental mathematically rigorous inquiries of the structure of solutions of the equations adopted as models in atmospheric dynamics have not yet been thoroughly done even for the simplest model. In this article we try to investigate the fundamental mathematical properties of the simplest model with the compressible Euler equations and the barotropic, or, isentropic motion of the ideal gas described in the co-ordinate system which rotates with a constant angular velocity. Although there remains open problems to be solved in mathematically rigorous way, some elementary aspects of the inquiry are exhibited in this article. \\

In the following in this section the equations and the boundary condition
to be considered are going to be described, and the definition of some concepts of solutions are going to be given. In  Section 2 stationary solutions with comactly  supported density and zero relative velocity  will be discussed. The limit of allowable
magnitude of the angular velocity of the rotation should be noted.  In Section 3 Lagrangian co-ordinate description of the equations for the perturbations near the stationary solutions will be discussed, and in Section 4 we shall show that the linearized wave equation for the perturbations allows an application of the Hille-Yosida theory of existence of solutions. Higher order regularities of the solutions remains in an open problem. In Section 5 the so called 'variational principle' will be formurated and its efficiency will be discussed. Although the concept of eigenfrequency and eigenvectors to the wave equations for the perturbations is remarkable, the existence and completeness of eigenvectors is in an open problem when the rotation is present. \\

Let us describe the situation to be considred precisely.

We consider motions of an atmosphere governed by the compressible Euler equations described by the uniformly rotating coordinate system $(t, \mbox{\boldmath$x$})$:
\begin{subequations}
\begin{align}
&\frac{D\rho}{Dt}+\rho(\nabla|\mbox{\boldmath$v$})=0, \label{1.1a}\\
&\rho\Big[\frac{D\mbox{\boldmath$v$}}{Dt}+2\mbox{\boldmath$\Omega$}\times\mbox{\boldmath$v$}
\Big]
+\nabla P+\rho\nabla \Phi^{\langle \Omega \rangle}=0 \label{1.1b}
\end{align}
\end{subequations}
on $t \in \mathbb{R}, \mbox{\boldmath$x$}=(x^1,x^2,x^3) \in 
\mathbb{R}^3\setminus {\mathfrak{B}}_0=
\{ \mbox{\boldmath$x$} \in \mathbb{R}^3 \  |\  r:=\|\mbox{\boldmath$x$}\| > R_0 \}$,
where
\begin{equation}
{\mathfrak{B}}_0=\{ \mbox{\boldmath$x$} \in \mathbb{R}^3 \  |\  r:=\|\mbox{\boldmath$x$}\| \leq R_0 \},
\end{equation}
$R_0$ being a positive number. The variables 
$\rho, P$ are the density, the pressure
and $\mbox{\boldmath$v$}=(v^1,v^2,v^3)^{\top}$ is the velocity field. We are denoting
\begin{equation}
 \frac{D}{Dt}:=
\frac{\partial}{\partial t}+(\mbox{\boldmath$v$}|\nabla)=
 \frac{\partial}{\partial t} +
\sum_{k=1}^3v^k\frac{\partial}{\partial x^k}
\end{equation}
and $\displaystyle (\nabla|\bm{v})=\sum_k
\frac{\partial v^k}{\partial x^k}$.
On the other hand
\begin{equation}
\mbox{\boldmath$\Omega$}=\Omega\frac{\partial}{\partial x^3},
\end{equation}
$\Omega$ being a constant, and $\Phi^{\langle \Omega \rangle}$ is the geopotential given as
\begin{equation}
\Phi^{\langle \Omega \rangle}(\mbox{\boldmath$x$})=-\frac{{\mathsf{G}M_0}}{r}-\frac{\Omega^2}{2}\varpi^2 \label{Geopote}
\end{equation}
where
$$ r=\sqrt{(x^1)^2+(x^2)^2+(x^3)^2},\quad \varpi=\sqrt{(x^1)^2+(x^2)^2},$$
${\mathsf{G}M_0}$ being a positive constant.
Since we are concerned with the value of the potential $\Phi^{\langle\Omega\rangle}$ only on the domain $\mathbb{R}^3\setminus \mathfrak{B}_0$, we assume that $\Phi^{\langle\Omega\rangle} \in C^{\infty}(\mathbb{R}^3)$ and \eqref{Geopote} holds for $r=\|\bm{x}\| >R_0$ by removing the singularity $\displaystyle \frac{1}{r}$ at $r=0$
by a smooth function near $r=0$, namely, we consider, instead of \eqref{Geopote}, 
\begin{equation}
\Phi^{\langle \Omega \rangle}(\bm{x})=
\chi\Big(\frac{r}{R_0}\Big)\Big[
-\frac{\mathsf{G}M_0}{r}-\frac{\Omega^2}{2}\varpi^2\Big],
\end{equation}
where $\chi \in C^{\infty}(\mathbb{R})$ such that
$1_{[\frac{1}{2}, +\infty[}
\leq \chi \leq 1_{[\frac{1}{4}, +\infty[}$.

The atmosphere is surrounding the Earth of radius $R_0$, mass $M_0$ and   $\mathsf{G}$ is the gravitational constant. $\Omega$ is the angular velocity of the rotation.\\

We assume \\

{\bf (A) :} {\it $P$ is the function of $\rho$ defined by
\begin{equation}
P=\mathsf{A}\rho^{\gamma} \quad\mbox{for}\quad \rho >0,
\end{equation}
$\mathsf{A}, \gamma$ being positive constants such that $1<\gamma<2$.
}\\

Under this assumption we introduce the variable $\varUpsilon$ by
\begin{equation}
\varUpsilon:=\int_0^{\rho}\frac{dP}{\rho}=\frac{\mathsf{A}\gamma}{\gamma-1}\rho^{\gamma-1}
\quad\mbox{for}\quad \rho >0.
\end{equation}\\

Since we are interested in  solutions with $\rho$ which has a vacuum region on which $\rho=0$, the concept of classical solutions needs a precise definition. The following are the definitions adopted in this article.\\

\begin{Definition}
An open connected subset $\mathfrak{D}$ of 
$\mathbb{R}^3$ is said to be an admissible domain with cover $\tilde{\mathfrak{D}}$  if $\tilde{\mathfrak{D}}=\mathfrak{D}\cup \mathfrak{B}_0$ is an open connected subset  of $\mathbb{R}^3$ such that 
$\mathfrak{D}=\tilde{\mathfrak{D}} \setminus \mathfrak{B}_0$ 
and $ \partial\mathfrak{D}$ is smooth.
\end{Definition}

Here and hereafter we use the following 
\begin{Notation}
Let $A, B$ be subsets of $\mathbb{R}^3$. $\complement A$ denotes the complement of $A$ and $A \setminus B$ stands for $A \cap \complement B$. $\mathbf{Cl}A$ stands for the closure of $A$, $\mathbf{Int}A$ stands for the interior of $A$, and $\partial A$ stands for the boundary of $A$, namely, $\partial A=\mathbf{Cl}A\setminus\mathbf{Int}A$.
\end{Notation}

Note that, if $\mathfrak{D}$ is an admissible domain, then, for a  sufficiently small positive $\delta$, we have 
$$\{ R_0 < r \leq (1+\delta)R_0 \} \subset \mathfrak{D}.$$

\begin{Definition}
Let $\mathfrak{D}$ be an admissible domain with cover $\tilde{\mathfrak{D}}$. A vector field $\bm{v}$ defined on $[0,T[\times \mathfrak{D}$ is said to be a classical admissible velocity field on $[0,T[\times\mathfrak{D}$ if  1) $\bm{v}$ has an extension onto $[0,T[\times\tilde{\mathfrak{D}}$ of class $C^1([0,T[\times \tilde{\mathfrak{D}})$
and 2) the boundary condition 
\begin{equation}
(\bm{n}|\bm{v})=0\quad\mbox{on}\quad [0,T[\times \partial \mathfrak{B}_0 \label{adminBC}
\end{equation}
holds, where $\mbox{\boldmath$n$}=
\mbox{\boldmath$n$}(\mbox{\boldmath$x$})=\mbox{\boldmath$x$}/R_0$ is the inward normal vector at the boundary point
$\mbox{\boldmath$x$} \in \partial\mathfrak{B}_0$.
\end{Definition}

\begin{Definition}
Let $\mathfrak{D}$ be an admissible doman with cover $\tilde{\mathfrak{D}}$. 
A function $(\rho, \mbox{\boldmath$v$})$ defined on $[0,T[\times \mathfrak{D}$ is said to be a classical
$(\rho, \mbox{\boldmath$v$})$-solution on $[0,T[\times \mathfrak{D}$ if
 1) $(\rho, \mbox{\boldmath$v$})$
has an extension onto $[0,T[\times \tilde{\mathfrak{D}}$ of class
$C^1([0,T[\times \tilde{\mathfrak{D}})$, 2) $\rho \geq 0$ on $[0,T[\times \mathfrak{D}$, 3) $(\rho, \mbox{\boldmath$v$})$ satsfies \eqref{1.1a}\eqref{1.1b} on $[0,T[\times \mathfrak{D}$, and 4) $\rho > 0$ on $[0,T[\times
\partial\mathfrak{B}_0$ and the boundary condition
\begin{equation}
(\bm{n}|\bm{v})=0 \quad\mbox{on}\quad [0,T[\times \partial\mathfrak{B}_0, \label{C.1}
\end{equation}
holds, where $\mbox{\boldmath$n$}=
\mbox{\boldmath$n$}(\mbox{\boldmath$x$})=\mbox{\boldmath$x$}/R_0$ is the inward normal vector at the boundary point
$\mbox{\boldmath$x$} \in \partial\mathfrak{B}_0$.
\end{Definition}

We want to permit the variable $\varUpsilon$ to take negative values somewhere. In order to do this we extend the equations devided by $\rho$ where $\rho >0$ to the equations:

\begin{subequations}
\begin{align}
&\frac{D\varUpsilon}{Dt}+(\gamma-1)\varUpsilon(\nabla|\mbox{\boldmath$v$})=0, \label{C.2a} \\
&\frac{D\mbox{\boldmath$v$}}{Dt}+2\mbox{\boldmath$\Omega$} \times \mbox{\boldmath$v$}
+\nabla(\varUpsilon +\Phi^{\langle\Omega\rangle})=0 \label{C.2b}
\end{align}
\end{subequations}\\

\begin{Definition}\label{DEF.3}
Let $\mathfrak{D}$ be an admissible domain with cover $\tilde{\mathfrak{D}}$. A function $(\varUpsilon, \mbox{\boldmath$v$})$
defined on $[0,T[\times\mathfrak{D}$ is said to be a classical
$(\varUpsilon,\mbox{\boldmath$v$})$-solution on $[0,T[\times\mathfrak{D}$, if
1) $(\varUpsilon, \mbox{\boldmath$v$})$
has an extension onto $[0,T[\times \tilde{\mathfrak{D}}$ of class
$C^1([0,T[\times \tilde{\mathfrak{D}})$, 
2)
$(\varUpsilon, \mbox{\boldmath$v$})$ satisfies \eqref{C.2a}\eqref{C.2b} on $[0,T[\times
\mathfrak{D}$, and 4) 
 $\varUpsilon >0$ on
$[0,T[\times \partial\mathfrak{B}_0$ and
the boundary condition
\begin{equation}
(\bm{n}|\bm{v})=0\quad\mbox{on}\quad [0,T[\times \partial\mathfrak{B}_0 \label{1.10}
\end{equation}
holds, where $\mbox{\boldmath$n$}=
\mbox{\boldmath$n$}(\mbox{\boldmath$x$})=\mbox{\boldmath$x$}/R_0$ is the inward normal vector at the boundary point
$\mbox{\boldmath$x$} \in \partial\mathfrak{B}_0$.
\end{Definition}

As for the variable $\varUpsilon$, $\varUpsilon \geq 0$ is not assumed. When 
$(\varUpsilon, \mbox{\boldmath$v$})$ is a classical $(\varUpsilon, \mbox{\boldmath$v$})$-solution on the admmisible domain $[0,T[\times \mathfrak{D}$, then $(\rho, \mbox{\boldmath$v$})$ defined by
\begin{equation}
\rho=\Big(\frac{\gamma-1}{\mathsf{A}\gamma}\varUpsilon\vee 0\Big)^{\frac{1}{\gamma-1}} \label{C.4}
\end{equation}
turns out to be a classical $(\rho, \mbox{\boldmath$v$})$-solution
on $[0,T[\times\mathfrak{D}$. Here and hereafter we use the

\begin{Notation}
We denote
\begin{equation}
Q\vee Q'=\max\{ Q, Q'\},\quad Q\wedge Q'=\min\{ Q, Q'\}.
\end{equation}
\end{Notation}

But the inverse may be impossible, namely, we cannot expect that any classical $(\rho, \mbox{\boldmath$v$})$-solution comes from classical $(\varUpsilon,\mbox{\boldmath$v$})$-solution by the above procedure, \eqref{C.4}, since
$1<\frac{1}{\gamma-1}<+\infty$ but $0<\gamma-1 <1$ so that
$\rho \in C^1$ does not imply $\rho^{\gamma-1}\in C^1$ at the vacumm boundary. 
In this sense the concept of classical $(\rho, \mbox{\boldmath$v$})$-solutions is weaker than that of classical $(\varUpsilon, \mbox{\boldmath$v$})$-solutions. 

Moreover we note that the uniqueness of solution cannot be expected under the concept of classical $(\rho, \mbox{\boldmath$v$})$-solutions in the following sense:

Let $(\rho, \mbox{\boldmath$v$})$ be a classical $(\rho, \mbox{\boldmath$v$})$-solution on $[0,T[\times
\mathfrak{D}$, $\mathfrak{D}$ being an admissible domain with  cover $\tilde{\mathfrak{D}}$. Suppose that $(]0,T[\times \mathfrak{D})
\setminus \mathbf{Cl}\{ \rho >0\}\not=\emptyset$, and consider a velocity field
$\mbox{\boldmath$v$}'\in C^1([0,T[\times \tilde{\mathfrak{D}})$ such that
$\mbox{\boldmath$v$}'(t,\mbox{\boldmath$x$})=\mbox{\boldmath$v$}(t,\mbox{\boldmath$x$})$ on $([0,T[\times \mathfrak{D})\setminus \mathcal{V}$ and
$\mbox{\boldmath$v$}'(t,\mbox{\boldmath$x$})\not=\mbox{\boldmath$v$}(t,\mbox{\boldmath$x$})$ on $\mathcal{V}$,
where $\mathcal{V}$ is a non-empty open subset of $(]0,T[\times\mathfrak{D})\setminus
\mathbf{Cl}\{\rho>0\}$. Clearly $(\rho, \mbox{\boldmath$v$}')$ is another classical
 $(\rho, \mbox{\boldmath$v$})$-solution such that it coincides with the original $(\rho, \mbox{\boldmath$v$})$ on $[0,t_1]\times\mathfrak{D}$ but it is different from the original $(\rho, \mbox{\boldmath$v$})$ after $t=t_1$ with $0<\exists t_1 <T$. \\


We are interested in motions with compactly supported $\rho$. Namely we use
\begin{Definition}
Let $\mathfrak{D}$ be an admisible domain. A classical $(\rho,\bm{v})$-solution $(\rho,\bm{v})$ on $[0,T[\times \mathfrak{D}$ is said to be compactly supported if,  for any fixed $t \in [0,T[$, $\mathbf{Cl}\{ \bm{x} | \rho(t,\bm{x})>0\}$
 is a comoact subset of $\partial\mathfrak{B}_0 \cup \mathfrak{D}$.
\end{Definition}

\section{Stationary solutions}

We are lokking for stationary solutions $\varUpsilon=\varUpsilon(\bm{x}), \bm{v}=\bm{v}(\bm{x})$. The equation to be satisfied are
\begin{subequations}
\begin{align}
(\bm{v}|\nabla)\varUpsilon +(\gamma-1)\varUpsilon(\nabla|\bm{v})=0, \label{601a} \\
(\bm{v}|\nabla)\bm{v}+2\bm{\Omega}\times\bm{v}+
\nabla(\varUpsilon+\Phi^{\langle\Omega\rangle})=0. \label{601b}
\end{align}
\end{subequations}
Let us use the cylindrical co-ordinate system $(\varpi, \phi, z)$ defined by
$$x^1=\varpi\cos\phi,\quad  x^2=\varpi\sin\phi, \quad  x^3=z.$$
The basis of the co-ordinates consists of the unit vectors
\begin{align*}
\bm{e}_{\varpi}&=\frac{1}{\varpi}\frac{\partial}{\partial\varpi}=
\frac{1}{\varpi}\Big(\cos\phi\frac{\partial}{\partial x^1}+\sin\phi\frac{\partial}{\partial x^2}\Big), \\
\bm{e}_{\phi}&=\frac{1}{\varpi}\frac{\partial}{\partial\phi}=-\sin\phi\frac{\partial}{\partial x^1}+\cos\phi\frac{\partial}{\partial x^2}, \\
\bm{e}_z&=\frac{\partial}{\partial x^3}.
\end{align*}\\

Suppose that the velocity field $\bm{v}$ is of th form
\begin{equation}
\bm{v}=V^{\phi}\bm{e}_{\phi}=\frac{V^{\phi}}{\varpi}\frac{\partial}{\partial\phi}.
\label{603}
\end{equation}
Since
$$(\bm{v}|\nabla)=\frac{V^{\phi}}{\varpi}\frac{\partial}{\partial\phi},\quad
(\nabla|\bm{v})=\frac{1}{\varpi}\frac{\partial V^{\phi}}{\partial\phi}, $$
the equation of continuity \eqref{601a} reduces to
$$\frac{V^{\phi}}{\varpi}\frac{\partial\varUpsilon}{\partial\phi}+
(\gamma-1)\frac{\varUpsilon}{\varpi}\frac{\partial V^{\phi}}{\partial\phi}=0. $$
This equation holds if 
\begin{equation}
\frac{\partial V^{\phi}}{\partial\phi}=0, \quad \frac{\partial\varUpsilon}{\partial\phi}=0. \label{604}
\end{equation}
So, supposing \eqref{604}, we solve te equation of motion \eqref{601b}. Since
\begin{align*}
(\bm{v}|\nabla)\bm{v}&=-(V^{\phi})^2\bm{e}_{\varpi}
+\frac{V^{\phi}}{\varpi}\frac{\partial V^{\phi}}{\partial\phi}\bm{e}_{\phi}=-(V^{\phi})^2\bm{e}_{\varpi}, \\
2\bm{\Omega}\times\bm{v}&=-2\Omega\varpi V^{\phi}\bm{e}_{\varpi}, \\
\nabla(\varUpsilon+\Phi^{\langle \Omega\rangle})&=
\frac{\partial}{\partial\varpi}(\varUpsilon+\Phi^{\langle \Omega\rangle})\bm{e}_{\varpi}
+\frac{1}{\varpi}\frac{\partial}{\partial\phi}(\varUpsilon+\Phi^{\langle \Omega\rangle})\bm{e}_{\phi} + \\
&+\frac{\partial}{\partial z}(\varUpsilon+\Phi^{\langle \Omega\rangle})
\bm{e}_z,
\end{align*}
the equation \eqref{601b} reduces to
\begin{subequations}
\begin{align}
-(V^{\phi})^2-2\Omega\varpi V^{\phi}+\frac{\partial}{\partial\varpi}(\varUpsilon+\Phi^{\langle \Omega\rangle})&=0, \label{605a} \\
\frac{\partial}{\partial z}(\varUpsilon+\Phi^{\langle \Omega\rangle})=0. \label{605b}
\end{align}
\end{subequations}
Here recall that $\frac{\partial}{\partial\phi}\Phi^{\langle\Omega\rangle}=0$ and that we are supposing $\frac{\partial}{\partial\phi}\varUpsilon=0$.\\

Taking
$$ \frac{\partial}{\partial z}\mbox{\eqref{605a}}-\frac{\partial}{\partial\varpi}\mbox{\eqref{605b}}=0,$$
we see $\partial V^{\phi}/\partial z=0$. Therefore there should exist a function $\omega$ such that
${V}^{\phi}=\varpi\omega(\varpi)$, namely, we consider
\begin{equation}
{\mbox{\boldmath$v$}}=\omega(\varpi)\frac{\partial}{\partial\phi}
=\varpi\omega(\varpi)\mbox{\boldmath$e$}_{\phi}.
\end{equation}

Integration of \eqref{605a}, \eqref{605b} gives
\begin{equation}
 \varUpsilon+
\Phi^{\langle \Omega \rangle}
+\frac{\Omega^2}{2}\varpi^2
-\mathsf{B}(\varpi)=\mbox{Const.}, \label{IntB}
\end{equation}
where
\begin{equation}
\mathsf{B}(\varpi):=\int_0^{\varpi}
(\omega(\acute{\varpi})+\Omega)^2\acute{\varpi}d\acute{\varpi}.
\end{equation}
Recall
$$\Phi^{\langle \Omega \rangle}
+\frac{\Omega^2}{2}\varpi^2=-\frac{{\mathsf{G}M_0}}{r}. $$
Therefore \eqref{IntB} reads 
\begin{equation}
\Upsilon + \Phi^{\langle \Omega, \omega\rangle}=\mbox{Const.},
\end{equation}
where
\begin{align}
\Phi^{\langle\Omega,\omega\rangle}&:=-\frac{\mathsf{G}M_0}{r}-\mathsf{B}(\varpi) \nonumber \\
&=-\frac{\mathsf{G}M_0}{r}
-\int_0^{\varpi}
(\omega(\acute{\varpi})+\Omega)^2\acute{\varpi}d\acute{\varpi}.
\end{align}

As in Section 2, let us specify the constant so that
\begin{equation}
{\varUpsilon}={\mathsf{G}M_0}\Big(\frac{1}{r}-\frac{1}{R}\Big)+\mathsf{B}(\varpi).
\end{equation}\\

Thus we have stationary soltion
\begin{equation}
\varUpsilon=\mathsf{G}M_0\Big(\frac{1}{r}-\frac{1}{R}\Big)+\mathsf{B}(\varpi), \qquad
\bm{v}=\varpi\omega(\varpi)\bm{e}_{\phi}.
\end{equation}
Here, in order to fix the idea, let us suppose the following assumption:\\

{\bf (B):} \  The function $
\omega $ belongs to the class $ C^1([0,+\infty[)\cap L^{\infty}(0,+\infty)$. \\


Now let us consider the static  stationary solution $(\rho, \bm{v})=
(\rho(\bm{x}), \mathbf{0})$: 
that is, $\omega(\varpi)=0$, and 
\begin{align}
\varUpsilon:&=
-\Phi^{\langle \Omega \rangle}-\frac{{\mathsf{G}M_0}}{R}  \\
\mbox{so that}& \nonumber \\
\varUpsilon&={\mathsf{G}M_0}\Big(\frac{1}{r}-\frac{1}{R}\Big)+\frac{\Omega^2}{2}\varpi^2
\quad\mbox{for}\quad r > R_0, \nonumber
\end{align}
where $R$ is a positive constant at the present. 
If $R \geq R_0$, then $R-R_0$ is the height of the stratosphere at the North and South Poles, and 
\begin{equation}
\varUpsilon_{\mathsf{P}}:=
{\mathsf{G}M_0}\Big(\frac{1}{R_0}-\frac{1}{R}\Big)\geq 0
\end{equation}
gives the density of the atmosphere at the Poles
\begin{equation}
\rho_{\mathsf{P}}:=\Big(\frac{\gamma-1}{\mathsf{A}\gamma}
\varUpsilon_{\mathsf{P}}\Big)^{\frac{1}{\gamma-1}}
=\Big(\frac{(\gamma-1){\mathsf{G}M_0}}{\mathsf{A}\gamma}\Big(\frac{1}{R_0}-\frac{1}{R}\Big)\Big)^{\frac{1}{\gamma-1}} .
\end{equation}\\

Now we are going to observe the shape of the set
\begin{equation}
\mathfrak{U}:=\{ (\varpi, z)\ |\  
0\leq\varpi, |z|<\infty, R_0 <r, 
-\Phi^{\langle \Omega \rangle}-\frac{{\mathsf{G}M_0}}{R}>0 \}.
\end{equation}

Let us introduce the non-dimensional variables $X=\varpi/R_0, Z=z/R_0$
and put
\begin{align}
&F(X^2, Z^2;\kappa):=\frac{1}{\sqrt{X^2+Z^2}}+\kappa X^2, \\
&\kappa=\frac{\Omega^2R_0^3}{2{\mathsf{G}M_0}}.
\end{align}
Then 
\begin{equation}
\Phi^{\langle \Omega \rangle}=-\frac{{\mathsf{G}M_0}}{R_0}F(X^2, Z^2; \kappa)
\end{equation}
and $\{ 
\displaystyle -\Phi^{\langle \Omega \rangle}-\frac{{\mathsf{G}M_0}}{R}
 > 0\}=\{ \displaystyle F >\frac{R_0}{R}\}$. Let us observe te shape of the level set $\{ F 
(X^2,Z^2;\kappa)=\lambda\}$ of $F$, $\lambda$ being a positive number. To do so, solving
\begin{equation}
F(X^2, Z^2; \kappa)=\lambda,
\end{equation}
we consider
\begin{equation}
Z^2=\frac{g(X^2)}{(\lambda-\kappa X^2)^2},
\end{equation}
for $\displaystyle  X^2 <{\frac{\lambda}{\kappa}}$, 
where
\begin{equation}
g(Q)=g(Q;\kappa, \lambda):=1-\lambda^2Q+2\lambda\kappa Q^2
-\kappa^2Q^3.
\end{equation}
We see 
$$g(0)=1,\qquad g\Big(\frac{\lambda}{\kappa}\Big)=1,$$
and
$$Dg(Q)=
(\lambda -\kappa Q)(3\kappa Q -\lambda).
$$
Note that $$g\Big(\frac{\lambda}{3\kappa}\Big)>0
\quad \Leftrightarrow \lambda^3 <\frac{27}{4}\kappa.$$\\

We see the shape of the graph of the function $g$ according to the following three cases:\\

Case(L): $\lambda^3 <\frac{27}{4}\kappa$.

Then $g(Q)>0$ for $0\leq Q <Q_{\infty}$, $g(Q_{\infty})=0$, 
$g(Q)<0$ for $Q_{\infty}<Q<+\infty$. Here $Q_{\infty}=Q_{\infty}(\frac{\lambda^3}{\kappa})$ is a number such that $Q_{\infty} > \frac{\lambda}{\kappa}$.\\

Case (M): $\lambda^3=\frac{27}{4}\kappa$.

Then $g(Q)>0$ for $0\leq Q <\frac{\lambda}{3\kappa}(=\frac{9}{4\lambda^2})$,
$g(\frac{\lambda}{3\kappa})=0$, $g(Q)>0$ for $\frac{\lambda}{3\kappa} <Q < Q_{\infty}$,
$g(Q_{\infty})=0$, $g(Q)<0$ for $Q_{\infty}<Q <+\infty$. \\

Case (H): $\lambda^3 >\frac{27}{4}\kappa$.

Then $g(Q)>0$ for $ 0\leq Q <Q_-$, $g(Q_-)=0$, $g(Q)<0$ for $Q_-<Q<Q_+$,
$g(Q_+)=0$, $g(Q)>0$ for $Q_+ < Q <Q_{\infty}$,
$g(Q_{\infty})=0$, $g(Q)<0$ for $Q_{\infty}<Q <+\infty$. Here $Q_{\pm}=Q_{\pm}(\frac{\lambda^3}{\kappa})$ are numbers such that
$0 < Q_- < \frac{\lambda}{3\kappa}<Q_+<\frac{\lambda}{\kappa}$ and
$Q_{\pm}\rightarrow \frac{\lambda}{3\kappa}$ as $\frac{\lambda^3}{\kappa}\rightarrow
\frac{27}{4}$.\\

Recall that we need the function $g(Q)$ only for $0\leq Q \leq \frac{\lambda}{\kappa}$, and we do not take care of the behavior of $g(Q)$
beyond $Q=\frac{\lambda}{\kappa}$, a fortiori, neither near nor beyond $Q=Q_{\infty}$. 
Anyway, correspondingly we see the shape of the set $\{  F >\lambda\}$ as follows:\\

In Case (L), the set $\{ F >\lambda\}$ is unbounded.\\

In Case (M), the set $\{ F > \lambda\}$ consistes of two connected components, say, $\mathfrak{U}_0$ and $\mathfrak{U}_{\infty}$, where
$\mathfrak{U}_0$ is bounded and included in $\{ 0 \leq X <\sqrt{\frac{\lambda}{3\kappa}}(=\frac{3}{2\lambda}) \}$ but $\mathfrak{U}_{\infty}$ is an unbounded subset of $\{ \sqrt{\frac{\lambda}{3\kappa}} < X< \sqrt{\frac{\lambda}{\kappa}} \}$.  Note that $Dg(\frac{\lambda}{3\kappa})=0$ so that $\partial\mathfrak{U}_0=\{ F=\lambda, X \leq \frac{3}{2\lambda} \}$ has a corner at the point $(X,Z)=(\sqrt{\frac{\lambda}{3\kappa}}, 0)$. 
Namely $\partial\mathfrak{U}_0$ near this point can be described as
$$Z=\pm C\Big(\frac{3}{2\lambda}-X\Big)
\Big(1+\Big[X-\frac{3}{2\lambda}\Big]_1\Big)\quad\mbox{as}\quad
X \rightarrow \frac{3}{2\lambda}-0,
$$
where $C=C(\kappa)$ is a positive constant and $[Y]_1$ stads for a convergent power series of the form
$\sum_{k\geq 1}a_kY^k$. \\

In Case (H), the set $\{ F >\lambda \}$ consistes of two connected components, say
$\mathfrak{U}_0 \subset \{ 0\leq X<\sqrt{Q_- }\}$ and
$\mathfrak{U}_{\infty} \subset
\{ \sqrt{Q_+} < X < \sqrt{\frac{\lambda}{\kappa}}  \}$. $\mathfrak{U}_0$ is bounded but $\mathfrak{U}_{\infty}$ is unbounded. Note that $X$ along $\partial\mathfrak{U}_0$ can be solved as a smooth ( real analytic) function of $Z, |Z| \ll 1$,
near the point $(X,Z)=(\sqrt{Q_-}, 0)$.
Namely the shape of $\partial\mathfrak{U}_0=\{ F=\lambda, X <\sqrt{\frac{\lambda}{3\kappa} } \}$ near this point is 
$$ Z=\pm C \sqrt{Q_--X}(1+[X-Q_-]_1)
\quad\mbox{as}\quad X \rightarrow Q_--0,$$
$C=C(\kappa,\lambda)$ beng a positive constant.
 The point is not a corner, and $\partial\mathfrak{U}_0$ is smooth. 
However it should be noticed that $\mathfrak{U}_0$ is not an ellipse.\\

Therefore, applying the above observations to $\displaystyle \lambda=\frac{R_0}{R}$, we claim

\begin{Theorem}
Compactly supported axially and equatorially symmetric static stationary solutions with the height of the stratosphere at the North Pole $R > R_0$ and constant angular velocity $\Omega$ exist if and only if  $\displaystyle\kappa= \frac{\Omega^2R_0^3}{2{\mathsf{G}M_0}} \leq \frac{4}{27}\Big(\frac{R_0}{R}\Big)^3$.
\end{Theorem} 

In this sense we should restrict ourselves to the case  $0 \leq \kappa < \frac{4}{27}\lambda^3$, that is, we require the following assumption\\

{\bf (K): } {\it $R > R_0$ and it holds 
\begin{equation}
0 \leq \Big(\frac{R}{R_0}\Big)^3\kappa=\frac{\Omega^2R^3}{2{\mathsf{G}M_0}} < \frac{4}{27}.
\end{equation}
}\\

We put
\begin{equation}
\mathcal{R}_1(\Omega^2):=
\begin{cases}
+\infty \quad\mbox{for}\quad \Omega^2=0 \\
\\
\displaystyle \frac{2}{3}\Big(\frac{{\mathsf{G}M_0}}{\Omega^2}\Big)^{\frac{1}{3}} \quad\mbox{for}\quad \Omega^2 > 0
\end{cases}.
\end{equation}
Then the condition {\bf (K)} reads
\begin{equation}
R_0 < R < \mathcal{R}_1(\Omega^2).
\end{equation}\\

In order that $R_0 <\mathcal{R}_1(\Omega^2)$ it is necessary that
\begin{equation}
\Omega^2 <\Omega^2_{\mathrm{max}}:=\frac{4}{9}\frac{{\mathsf{G}M_0}}{(R_0)^3}.
\end{equation}\\


Under the assumption {\bf (K)}, or, in Case (H) with $\lambda=
\frac{R_0}{R}$, we
consider the number $\frac{3R}{2R_0}=\frac{3}{2\lambda}$. Since 
$\frac{3}{2\lambda}<\sqrt{ \frac{\lambda}{3\kappa}}$ and since  $\displaystyle g\Big(\frac{9}{4\lambda^2}\Big)<0$, which is not obvious but can be shown, we see that $\sqrt{Q_-} <\frac{3}{2\lambda}$. On the other hand, we see obviosly that 
$ \frac{3}{2\lambda}<\sqrt{\frac{\lambda}{3\kappa}}<\sqrt{Q_+}$. Therefore we have the estimate
$$ \sqrt{Q_-} < \frac{3}{2\lambda}=\frac{3R}{2R_0} <\sqrt{Q_+}. $$ Thus
the closure of  $\mathfrak{U}_0$ is included in
$\{R_0 \leq r,\varpi <{3R}/{2}\}$, and $\varUpsilon$ is a classical solution
on $[0,+\infty[\times \mathfrak{D}$, where
$\mathfrak{D}=\{R_0 < r,\varpi <{3R}/{2}\}$. So
$$\bar{\rho}=
\Big(\frac{\gamma-1}{\mathsf{A}\gamma}(
\varUpsilon \restriction \mathfrak{D} )
\vee 0\Big)^{\frac{1}{\gamma-1}}$$
turns out to be a classical solution on $[0,+\infty[\times\mathfrak{D}$. In other words, if we  put
\begin{equation}
\varUpsilon_{\blacktriangle}=
\begin{cases}
\varUpsilon \quad\mbox{on}\quad \mathfrak{U}_0 \\
0 \quad\mbox{on}\quad \{R_0 < r\}\setminus \mathfrak{U}_0,
\end{cases}
\end{equation}
that is,
\begin{align}
\varUpsilon_{\blacktriangle}&=(\varUpsilon\vee  0)\cdot 1_{\varpi <3R/2}
=(\varUpsilon\cdot 1_{\varpi < 3R/2})\vee 0 \nonumber \\
&=
\begin{cases}
\varUpsilon \quad\mbox{if}\quad \varUpsilon>0\quad\mbox{and}\quad \varpi <\frac{3R}{2} \\
0 \quad\mbox{otherwise}
\end{cases},
\end{align}
the density distribution $\bar{\rho}$  given by
\begin{equation}
\bar{\rho}=\Big(\frac{\gamma-1}{\mathsf{A}\gamma}
\varUpsilon_{\blacktriangle}\Big)^{\frac{1}{\gamma-1}}
\end{equation}
gives a classical solution $(\bar{\rho}, \mathbf{0})$ on $[0,\infty[\times
 ( \mathbb{R}^3\setminus {\mathfrak{B}}_0)$.
Note that $(\varUpsilon \cdot 1_{\varpi <3R/2}, \mathbf{0})$ is not a classical $(\varUpsilon, \mbox{\boldmath$v$})$-solution on $[0,+\infty[\times ( \mathbb{R}^3\setminus {\mathfrak{B}}_0)$ in the sense of Definition \ref{DEF.3}.\\

We are considering the uniformly rotating atmosphere  which occupies
\begin{equation}
\mathfrak{R}:=\{ \mbox{\boldmath$x$} \in \mathbb{R}^3 \  |\  
R_0 <r, \rho^{\flat}(\mbox{\boldmath$x$})=\rho(\varpi, z)>0 \},
\end{equation}
provided {\bf (K)}. 
Then  the boundary $\partial\mathfrak{R}$ of he atmosphere consists of the two connected components $\Sigma_0=\partial\mathfrak{B}_0=\{ r = R_0\}$, the surface of the Earth,
and $$\Sigma_1=\{ \frac{1}{\sqrt{X^2+Z^2}}+\kappa X^2=\frac{R_0}{R}, \quad\mbox{with}\quad
\varpi=R_0X, z=R_0Z \}, $$
the stratosphere of the atmosphere. We are going to
show that the boundary $\Sigma_1$ is a physical vacuum boundary, that is,
at each boundary point $\mathrm{P} \in \Sigma_1$ we have ($\nabla \varUpsilon^{\flat}|\mbox{\boldmath$n$})<0$, where $\mbox{\boldmath$n$}$ is the outer normal vector at $\mathrm{P}$ and $\varUpsilon^{\flat}(\mbox{\boldmath$x$})=\varUpsilon(\varpi, z)$. Recall that $\displaystyle \varUpsilon=\frac{1}{\gamma-1}\frac{dP}{d\rho}$, while $\displaystyle \frac{dP}{d\rho}$ is the square of the sound speed.

In fact if we consider the situation in the $(X, Z)$-plane at $\mathrm{P}: X=X_1, Z= Z_1$ with $0<X_1<\sqrt{Q_-}, Z_1=f(X_1)$ where
$$f(X):=\frac{\sqrt{g(X^2)}}{\lambda-\kappa X^2} \quad\mbox{with}\quad \lambda=\frac{R_0}{R}. $$ 
By a straight calculation, it can be shown that $Df(X) < 0$ for $0<X<\sqrt{Q_-}$. We have
$$\mbox{\boldmath$n$}=\frac{1}{1+Df(X_1)^2}\Big(-Df(X_1)\frac{\partial}{\partial X}+\frac{\partial}{\partial Z}\Big)$$ so that 
\begin{align*}
(\nabla \varUpsilon^{\flat}|\mbox{\boldmath$n$})&=
\frac{1}{1+Df(X_1)^2}\Big(-\frac{\partial \varUpsilon}{\partial X}Df(X_1)+\frac{\partial \varUpsilon}{\partial Z}\Big) \\
&=\frac{1}{1+Df(X_1)^2}\frac{{\mathsf{G}M_0}}{R_0}\Big(
\frac{XDf(X)-f (X)}{(X^2+f(X)^2)^{3/2}}-2\kappa XDf(X)\Big)_{X=X_1}.
\end{align*}
Since $0< X_1, 0<f(X_1), Df(X_1) <0, 0 <\kappa $, we see $(\nabla \varUpsilon^{\flat}|\mbox{\boldmath$n$}) <0$ at $\mathrm{P}$. The exceptional cases, the North Pole ($X=0$) and the Equator 
($X=\sqrt{Q_-}$\  ) , can be checked easily. \\

Now we consider a not static stationary solution
\begin{equation}
\varUpsilon=\mathsf{G}M_0\Big(\frac{1}{r}-\frac{1}{R}\Big)+\mathsf{B}(\varpi), \qquad
\bm{v}=\varpi\omega(\varpi)\bm{e}_{\phi}
\end{equation}
under the assumption:\\

{\bf (B):} \  The function $
\omega $ belongs to the class $ C^1([0,+\infty[)\cap L^{\infty}(0,+\infty)$. \\

Under this assumption, we put
\begin{equation}
\tilde{\kappa}(X^2)=
\frac{R_0}{\mathsf{G}M_0}
\frac{\mathsf{B}(R_0X)}{X^2}=
\frac{R_0}{\mathsf{G}M_0}\frac{1}{X^2}
\int_0^{R_0X}(\omega(\acute{\varpi})+\Omega)^2\acute{\varpi}d\acute\varpi.
\end{equation}
Note $X^2 \mapsto \tilde{\kappa}(X^2)X^2$ is continuos and montone
nondecreasing. 

Put
\begin{equation}
\kappa:=\sup_{X^2>0}\tilde{\kappa} =
\frac{R_0^3}{2\mathsf{G}M_0}\|\omega+\Omega\|^2,
\end{equation}
where
\begin{equation}
\|\omega+\Omega\|=\sup_{\varpi>0}|\omega(\varpi)+\Omega|.
\label{0408}
\end{equation}

Let us consider the stationary solution with compactly  supported $\rho$ under the assumption \\

($\tilde{\mathbf{K}}$): It holds that
\begin{equation}
R_0 <R,\qquad \Big(\frac{R}{R_0}\Big)^3\kappa
=\frac{R^3}{2\mathsf{G}M_0}\|\omega+\Omega\|^2<\frac{4}{27}.
\end{equation}
The density is given by
\begin{equation}
{\rho}=\Big(\frac{\gamma-1}{\mathsf{A}\gamma}\Big)^{\frac{1}{\gamma-1}}
({\varUpsilon}_{\blacktriangle})^{\frac{1}{\gamma1-1}}
\quad\mbox{with}\quad
{\varUpsilon}_{\blacktriangle}=({\varUpsilon}\vee 0)\cdot 1_{\varpi <3R/2}.
\end{equation}

In fact, $\{ \varUpsilon >0\}$ is $\{F(X^2,Z^2;\tilde{\kappa}) >\lambda=
\frac{R_0}{R} \}$, where
\begin{equation}
F(X^2,Z^2;\tilde{\kappa})=
\frac{1}{\sqrt{X^2+Z^2}}+\tilde{\kappa}(X^2)X^2.
\end{equation}
Here we are using the change of varialble 
$\varpi=R_0X, z=R_0Z$

Then $F(X^2, Z^2;\tilde{\kappa}) >\lambda$ if and only if
either
1) $\tilde{\kappa}(X^2)X^2 \geq \lambda$, or
2)
$$\tilde{\kappa}(X^2)X^2 <\lambda \quad \mbox{and}\quad Z^2<f(X^2;\tilde{\kappa})
$$
where 
$$f(X^2;\tilde{\kappa})=\frac{1}{(\lambda-\tilde{\kappa}(X^2)X^2)^2}
-X^2.
$$
Note that $f(X^2;\tilde{\kappa})\leq f(X^2;\kappa)$. 
Now we are supposing $$\tilde{\kappa}(X^2) \leq \kappa<
\frac{4}{27}\Big(\frac{R_0}{R}\Big)^3.$$

So, $\{ \varUpsilon >0\}$ has a bounded connected component of the form
$\{ 0<X^2 < Q_-, Z^2<f(X^2,\tilde{\kappa})\} $. Here 
$Q_-$ is a positive number  $\leq Q_-
(\frac{\lambda^3}{\kappa})<\frac{3}{2}\lambda$ such that
$f(X^2;\tilde{\kappa})>0$ for $0<X^2<Q_-$, and $f(Q_-;\tilde{\kappa})=0$.

\begin{Remark} 
In this situation with $R (>R_0)$ being fixed, it is sufficient that the function $\omega$ is given on the finite interval  $[0, {3R}/{2}]$ as a function of $C^1$-class on this interval, under a weaker assumption than {\bf (B)}, and we can put
$$
\|\omega+\Omega\|_{\infty}=\sup_{0\leq \varpi \leq 3R/2}|\omega(\varpi)+\Omega|
$$
instead of \eqref{0408}.
\end{Remark}

\section{Description of the flow by the Lagrangian co-ordinate}

Let $\mathfrak{D}_0$ be an admissible domain with cover $\tilde{\mathfrak{D}_0}$, and let $ \mbox{\boldmath$v$}$ be a classical admissible velocity field on $[0,T_0[\times \mathfrak{D}_0$.
Let $\mathfrak{D}$ with a cover $\tilde{\mathfrak{D}}$ be an admissible subdomain of $\mathfrak{D}$,
which we shall call a {\it proper asdmissible subdomain of} $\mathfrak{D}$, 
 such that 
 $\tilde{ \mathfrak{D}} \Subset \tilde{\mathfrak{D}_0}$, namely, 
 there is a compact set $K$ such that
$\mathbf{Cl}\tilde{\mathfrak{D}} \subset \mathbf{Int}K \subset K \subset \tilde{\mathfrak{D}_0}$.

Then {\it there is a positive number $T (<T_0)$ such that
for any $(\tau, \bar{\bm{x}}) \in [0,T]\times \tilde{\mathfrak{D}}$ the initial value problem of the ordinary differential equation
\begin{subequations}
\begin{align}
&\frac{d \bm{x}}{dt}=\mbox{\boldmath$v$}(t, \bm{x}), \label{FlowIVa}\\
& \bm{x}|_{t=\tau}=\bar{{\bm{x}}} \label{FlowIVb}
\end{align}
\end{subequations}
admits the unique solution $\bm{x}=\bm{\varphi}(t, \tau, \bar{{\bm{x}}})$ which exists on  $t \in [0, T]$ while $ \bm{x} \in \tilde{\mathfrak{D}_0}$}. We shall call $\varphi$ {\it the flow associated with the velocity field} $\bm{v}$. The flow $\varphi$, as a function, belongs to the class 
$C^1([0,T]\times[0,T]\times \mathbf{Cl}\tilde{\mathfrak{D}}\  ;\  \tilde{\mathfrak{D}_0})$.\\

In fact the fundamental theorems of ordinary differential equations (see e.g., \cite{CoddingtonL} ) applied to the vector field $\mbox{\boldmath$v$}$ considered as a $C^1$-field on $[0,T_0[\times \tilde{\mathfrak{D}_0}$ guarantee the existence and uniqueness of the flow on $[0,T]$
valued in $\tilde{\mathfrak{D}_0}$ for $\mbox{\boldmath$x$} \in \mathbf{Int}K$. On the other hand the boundary condition \eqref{adminBC} of $\mbox{\boldmath$v$}$ on $\partial \mathfrak{B}_0$ guarantees that $\partial \mathfrak{B}_0$ is an invarinat set of the equation so that the flow starting with
$\bar{\bm{x}} \in \mathfrak{D}$ cannot touch $\partial\mathfrak{B}_0$ and must remain inside of $\complement \mathfrak{B}_0$, namely,
$\bm{x} \in \tilde{\mathfrak{D}}\setminus \mathfrak{B}_0=\mathfrak{D}$.\\

Of course, abstractly speaking, we can consider the existence domain for $\bm{\varphi}$ of the form
$\mathcal{O}=\bigcup_{(\tau, \bar{\bm{x}}) \in [0,T_0[\times \tilde{\mathfrak{D}_0}} 
I_{\tau, \bar{\bm{x}}} \times\{\tau\}\times\{\bar{\bm{x}}\}$, where $I_{\tau, \bar{\bm{x}}}
\subset [0,T_0[$ is the maximal interval of existence of the solution of \eqref{FlowIVa}\eqref{FlowIVb}.
The above observation says that
$[0,T] \subset I_{\tau, \bar{\bm{x}}}$ uniformly for $\forall (\tau, \bar{\bm{x}} )\in [0,T]\times \tilde{\mathfrak{D}}$.\\

 Hereafter, we consider $\bm{\varphi}(t, 0, \bar{\bm{x}})$ and denote,
for the sake of brevity, 
 $$ \bm{\varphi}(t,\bar{\bm{x}}):=\bm{\varphi}(t,0, \bar{\bm{x}}).$$
 Thus 
 \begin{subequations}
 \begin{align}
 &\frac{\partial}{\partial t}\bm{\varphi}(t,\bar{\bm{x}})=\bm{v}(t,\bm{\varphi}(t,\bar{\bm{x}})), \label{3.2a}\\
 &\bm{\varphi}(0,\bar{\bm{x}})=\bar{\bm{x}},\label{3.2b}
 \end{align}
 \end{subequations}
 for $ (t, \bar{\bm{x}}) \in [0,T]\times \tilde{\mathfrak{D}}$. \\

Now let us suppose that the velocity field $\bm{v}$ is that of a classical $(\varUpsilon,\bm{v})$-solution $(\varUpsilon,\bm{v})$ on $[0,T_0[\times \mathfrak{D}_0$.
We consider the associated flow
$\bm{\varphi} \in C^1([0,T]\times \mathbf{Cl}\tilde{\mathfrak{D}})$ 
which satisfies \eqref{3.2a}\eqref{3.2b} for $\forall (t,\bar{\bm{x}}) \in [0,T]\times \tilde{\mathfrak{D}}$.\\

Hereafter we denote
\begin{equation}
(D_{\bar{\bm{x}}}{\bm{x}})(t, \bar{\bm{x}})=D\bm{\varphi} (t,\bar{\bm{x}})=\Big(\frac{\partial}{\partial \bar{\bm{x}}^{\alpha}}\Big(\bm{\varphi}(t,\bar{\bm{x}})\Big)^j\Big)_{j,\alpha}.
\end{equation}

We shall use the inverse matrix of $D_{\bar{\bm{x}}}\bm{x}(t,\bar{\bm{x}})$. Actually thereis a positive number 
$T_1(\leq T)$ such that 
$D_{\bar{\bm{x}}}\bm{x}(t,\bar{\bm{x}})=D\bm{\varphi}(t,\bar{\bm{x}})$ is invertible when $(t, \bar{\bm{x}}) \in [0,T_1]\times \tilde{\mathfrak{D}}$. ( Proof. Let
$$\Delta=D\bm{\varphi}(t,\bar{\bm{x}})-I=D\bm{\varphi}(t,\bar{\bm{x}})-D \bm{\varphi} (t,\bar{\bm{x}}).$$
Since $\bm{\varphi} \in C^1([0,T]\times \mathbf{Cl}\tilde{\mathfrak{D}})$, for any $\Theta , 0<\Theta<1,$ there is $T_1 \leq T$ such that $\|\Delta\|<\Theta$ for
$(t,\bar{\bm{x}}) \in [0,T_1]\times \tilde{\mathfrak{D}}$. then $(D\bm{\varphi}(t,\bar{\bm{x}}))^{-1}$ is
given by the Neumann series $\sum_{k=0}^{\infty}(-\Delta)^k$. $\square$. )

If we denote by $\bm{\psi}_t$ the inverse mapping of $\bm{\varphi}(t,\cdot)\restriction _{\mathfrak{V}}$, $\mathfrak{V}$ being a neighborhood of
$\bar{\bm{x}}_0 \in \tilde{\mathfrak{D}}, t \in [0,T_1]$, then for $\bar{\bm{x}} \in \mathfrak{V}$, $(D_{\bar{\bm{x}}}\bm{x}(t,\bar{\bm{x}}))^{-1}$ is nothing but  the matrix
 $\displaystyle \Big(\frac{\partial}{\partial x^j}\Big(\psi_t(\bm{x})\Big)^{\alpha}\Big)_{\alpha,j}$ at $\bm{x}=\bm{\varphi}(t,\bar{\bm{x}})$. \\


We note that it holds
\begin{equation}
D_{\bar{\bm{x}}}\bm{x}(t,\bar{\bm{x}})=D\bm{\varphi}(t,\bar{\bm{x}})=
\exp\Big[
\int_0^t(D_{\bm{x}}\bm{v})(s,\bm{\varphi}(s,\bar{\bm{x}}))ds\Big] \label{3.6}
\end{equation}
for $(t,\bar{\bm{x}}) \in [0,T]\times \tilde{\mathfrak{D}}$.  Here, of course,
$D_{\bm{x}}\bm{v}$ stands for the matrix
$\displaystyle \Big(\frac{\partial v^k}{\partial x^j}\Big)_{k,j}$. 

In fact, we see
\begin{align*}
\frac{\partial}{\partial t}D_{\bar{\bm{x}}}x(t,\bar{\bm{x}})&=\frac{\partial}{\partial t}D\bm{\varphi}(t,\bar{\bm{x}})= \\
&=D_{\bar{\bm{x}}}\frac{\partial}{\partial t}\bm{\varphi}(t,\bar{\bm{x}})=D_{\bar{\bm{x}}}\bm{v}(t,\bm{\varphi}(t,\bm{x}))= \\
&=D_{\bm{x}}\bm{v}(t,\bm{x})\Big|_{\bm{x}=\bm{\varphi}(t,\bar{\bm{x}})}. D_{\bar{\bm{x}}}\bm{x}(t,\bar{\bm{x}}).
\end{align*}
and $$
D_{\bar{\bm{x}}}(0, \bar{\bm{x}})=I.
$$.\\


Let us put
\begin{equation}
\overset{\circ}{\varUpsilon}(\bm{x})=\varUpsilon(0,\bm{x}).
\end{equation}
and
\begin{equation}
\varUpsilon^L(t,\bar{\bm{x}})=\varUpsilon(t,\bm{\varphi}(t,\bar{\bm{x}}))
\end{equation}
for $(t,\bar{\bm{x}})\in [0,T]\times \tilde{\mathfrak{D}}$.\\

  Integrating the equation of continuity \eqref{C.2a}, we claim
  
  \begin{Proposition}\label{Prop.Conti}
It holds that
\begin{equation}
\varUpsilon^L(t,\bar{\bm{x}})=
\varUpsilon(t, \bm{\varphi}(t,\bar{\bm{x}}))=
\overset{\circ}{\varUpsilon}(\bar{\bm{x}})\mathrm{det}
D\bm{\varphi}(t,\bar{\bm{x}})^{-(\gamma-1)}  \label{ECL}
\end{equation}
for $t \in [0,T], \bar{\bm{x}} \in \mathfrak{D}$.
\end{Proposition}



Proof. 
The equation \eqref{C.2a} reads
$$\frac{D}{Dt}\log \varUpsilon=-(\gamma-1)(\nabla_{\bar{x}}|\mbox{\boldmath$v$}) (t,{\mbox{\boldmath$x$}})$$
where $\varUpsilon \not=0$. Therefore
$$ \varUpsilon(t,\bm{x})=\overset{\circ}{\varUpsilon}(\bar{\bm{x}})\exp\Big[
-(\gamma-1)\int_0^t
(\nabla_{\bar{x}}|\mbox{\boldmath$v$})(t', \bm{\varphi}(t',\bar{\mbox{\boldmath$x$}}))dt' \Big]
$$
We note
$$
(\nabla_{\bar{x}}|\mbox{\boldmath$v$})(t,{\mbox{\boldmath$x$}})=\mathrm{tr}D_{\bar{x}}\mbox{\boldmath$v$}(t,{\mbox{\boldmath$x$}}).
$$
But, since
\begin{align*}
\Big(\frac{\partial}{\partial t}\Big)_{\bar{\mbox{\boldmath$x$}}}D\bm{\varphi}(t,\bar{\mbox{\boldmath$x$}})&=D_{\bar{\mbox{\boldmath$x$}}}\Big(\frac{\partial}{\partial t}\Big)_{\bar{\mbox{\boldmath$x$}}}{\mbox{\boldmath$x$}}(t,0, \bar{\mbox{\boldmath$x$}})=D_{\bar{\mbox{\boldmath$x$}}}\mbox{\boldmath$v$}(t,{\mbox{\boldmath$x$}}) \\
&=D_{{\mbox{\boldmath$x$}}}\mbox{\boldmath$v$}(t,{\mbox{\boldmath$x$}}).D_{\bar{\mbox{\boldmath$x$}}}{\mbox{\boldmath$x$}} =
D_{{\mbox{\boldmath$x$}}}\mbox{\boldmath$v$}(t,{\mbox{\boldmath$x$}}).D\bm{\varphi}(t,\bar{\mbox{\boldmath$x$}}),
\end{align*}
we have
$$D_{\bar{x}}\mbox{\boldmath$v$}(t,{\mbox{\boldmath$x$}})=
\Big(\frac{\partial}{\partial t}\Big)_{\bar{\mbox{\boldmath$x$}}}D\bm{\varphi}(t,\bar{\mbox{\boldmath$x$}}).
D\bm{\varphi}(t,\bar{\mbox{\boldmath$x$}})^{-1}.
$$
Thus
$$(\nabla_{\bar{x}}|\mbox{\boldmath$v$})(t,{\mbox{\boldmath$x$}})=
\mathrm{tr}\Big(\Big(\frac{\partial}{\partial t}\Big)_{\bar{\mbox{\boldmath$x$}}}D\bm{\varphi}(t,\bar{\mbox{\boldmath$x$}}).
D\bm{\varphi}(t,\bar{\mbox{\boldmath$x$}})^{-1} \Big)
=\Big(\frac{\partial}{\partial t}\Big)_{\bar{\bm{x}}}\log \mathrm{det} D\bm{\varphi}(t,\bar{\mbox{\boldmath$x$}}).
$$
Since $D\bm{\varphi}(0,\bar{\mbox{\boldmath$x$}})=I$, it follows that 
$$\int_0^t(\nabla_{\bar{x}}|\bm{v})(t', \varphi(t', \bar{\bm{x}}))dt'=
\log\mathrm{det}D\bm{\varphi}(t,\bar{\mbox{\boldmath$x$}}).$$
$\square$\\

Let us consider the equation of motion \eqref{C.2b}, namely
\begin{equation}
 \frac{D\bm{v}}{D t}+B\bm{v}+
\mathrm{glad}_{\bm{x}} (\varUpsilon+\Phi^{\langle\Omega\rangle})=0, \label{EM1}
\end{equation}
where  we denote
\begin{equation}
B\bm{v}=2\Omega 
\begin{bmatrix}
0 \\
\\
0 \\
\\
1\end{bmatrix}
\times \bm{v}=
2\Omega
\begin{bmatrix}
-v^2 \\
\\
v^1 \\
\\
0
\end{bmatrix}
=
2\Omega
\begin{bmatrix}
0 & -1 & 0 \\
\\
1 & 0 & 0 \\
\\
0 & 0 & 0
\end{bmatrix}
\bm{v}.
\end{equation}\\

We put
\begin{align}
\bm{v}^L &:=
\Big(D\bm{\varphi}(t,\bar{\bm{x}})\Big)^{-1}\bm{v}(t,\bm{\varphi}(t,\bar{\bm{x}})) \nonumber \\
&=
\Big(D\bm{\varphi}(t,\bar{\bm{x}})\Big)^{-1}\frac{\partial}{\partial t}\bm{\varphi}(t,\bar{\bm{x}}).
\end{align}

Then \eqref{EM1} reads
\begin{align}
&\frac{\partial \bm{v}^l}{\partial t}+(D_{\bar{\bm{x}}}\bm{x})^{-1}\frac{\partial}{\partial t}(D_{\bar{\bm{x}}}\bm{x}) \bm{v}^L + \nonumber \\
&+(D_{\bar{\bm{x}}}\bm{x})^{-1}B(D_{\bar{\bm{x}}}\bm{x})\bm{v}^L+
(D_{\bar{\bm{x}}}\bm{x})^{-1}
((D_{\bar{\bm{x}}}\bm{x})^{-1})^{\top}\mathrm{grad}_{\bar{\bm{x}}}
(\Upsilon^L+\Phi^{\langle\Omega\rangle L})=0 \label{EM2}
\end{align}
on $(t,\bar{\bm{x}}) \in [0,T]\times \tilde{\mathfrak{D}}$. Here
\begin{equation}
\Phi^{\langle\Omega\rangle L}(t, \bar{\bm{x}})
:=\Phi^{\langle\Omega\rangle }(\bm{\varphi}(t,\bar{\bm{x}})).
\end{equation}\\

Let us consider tge boundary condition
\begin{equation}
(\bm{n}(\bm{x})|\bm{v}(t,\bm{x}))=0\quad \forall (t,\bm{x}) \in
[0,T]\times \partial\mathfrak{B}_0, \label{BCE}
\end{equation}
where
\begin{equation}
\bm{n}(\bm{x})=\frac{1}{R_0}\bm{x}\quad\mbox{for}\quad 
\bm{x} \in \partial\mathfrak{B}_0=\{ \|\bm{x}\|=R_0\}.
\end{equation}

Put
\begin{equation}
\bm{n}^L(t,\bar{\bm{x}})=\frac{1}{R_0}(D_{\bar{\bm{x}}}\bm{x}(t,\bar{\bm{x}}))^{\top}\bm{\varphi}(t,\bar{\bm{x}})
\end{equation}
for $(t,\bar{\bm{x}}) \in [0,T]\times \partial\mathfrak{B}_0$. Then 
\eqref{BCE} reads
\begin{equation}
(\bm{n}^L(t,\bar{\bm{x}})|\bm{v}^L(t,\bar{\bm{x}}))=0
\quad\mbox{for}\quad
(t,\bar{\bm{x}})\in [0,T]\times \partial\mathfrak{B}_0. \label{BCL}
\end{equation}

Note that \eqref{BCL} means
\begin{equation}
\frac{1}{2}\frac{\partial}{\partial t}\|\bm{\varphi}(t,\bar{\bm{x}})\|^2
=\Big(\bm{\varphi} \Big|\frac{\partial}{\partial t}\bm{\varphi} \Big)=0 \quad\mbox{if}\quad \bar{\bm{x}} \in \partial\mathfrak{B}_0,
\end{equation}
so that , if $\bar{\bm{x}} \in \partial\mathfrak{B}_0$, then 
$\bm{\varphi}(t,\bar{\bm{x}}) \in \partial\mathfrak{B}_0$ for $\forall t $. \\


Now let us suppose that there is a vector field ${\bm{\varphi}}(t, \bar{\bm{x}}), (t, \bar{\bm{x}}) \in 
[0,T]\times \tilde{\mathfrak{D}} ,$ such that $ \bm{\varphi}
\in \bigcap_{\ell =1,2}C^{2-\ell}([0, T];C^{\ell}(\mathbf{Cl}\tilde{\mathfrak{D}} ))$ and
$$ \bm{\varphi}(0,\bar{\bm{x}})=\bar{\bm{x}} $$
which satisfies the equation \eqref{EM2},
where we read
\begin{subequations}
\begin{align}
&D_{\bar{\bm{x}}}\bm{x}=D\bm{\varphi}(t,\bar{\bm{x}}), \\
&\bm{v}^L=\Big(D\bm{\varphi}(t,\bar{\bm{x}})\Big)^{-1}\frac{\partial}{\partial t}\bm{\varphi}(t,\bar{\bm{x}}), \\
&\varUpsilon^L=\overset{\circ}{\varUpsilon}(\bar{\bm{x}})
\det \Big(D\bm{\varphi}(t,\bar{\bm{x}})\Big)^{-(\gamma-1)}, \\
&\Phi^{\langle \Omega \rangle L}=\Phi^{\langle \Omega \rangle}(\bm{\varphi}(t,\bar{\bm{x}})),
\end{align}
\end{subequations}
and
the boundary condition
\eqref{BCL}, where we read
\begin{equation}
\bm{n}^L(t,\bar{\bm{x}})=\frac{1}{R_0}(D\bm{\varphi}(t,\bar{\bm{x}}))^{\top}\bm{\varphi}(t,\bar{\bm{x}}),
\end{equation}
forgetting that the vector field $\bm{\varphi}$ comes from the flow generated by the velocity field $\bm{v}$ of a classical solution $(\varUpsilon,\bm{v})$.
Namely, we use the following definitions:

\begin{Definition}
Let $\mathfrak{D}$ be an admissible domain with cover $\tilde{\mathfrak{D}}$. A scalar field $\overset{\circ}{\varUpsilon}$ defined on $\tilde{\mathfrak{D}}$ is called an admissible $\varUpsilon$-data, if $\overset{\circ}{\varUpsilon} \in C^1(\mathbf{Cl}\tilde{\mathfrak{D}})$ and
$\overset{\circ}{\varUpsilon} >0$ on $\partial\mathfrak{B}_0$.
\end{Definition}
\begin{Definition}
Let $\mathfrak{D}$ be an admissible domain with cover $\tilde{\mathfrak{D}}$, and let $\overset{\circ}{\varUpsilon}$ be an admissible $\varUpsilon$-data on $\mathfrak{D}$. A vector field $\varphi$ is said to be an admissible flow on $[0,T]\times\mathfrak{D}$
associated with $\overset{\circ}{\varUpsilon}$, if 1) $ \bm{\varphi}
\in \bigcap_{\ell =1,2}C^{2-\ell}([0, T];C^{\ell}(\mathbf{Cl}\tilde{\mathfrak{D}} ))$, 
2) the equation \eqref{EM2} and the boundary condition \eqref{BCL} are satisfied on $[0,T]\times \mathfrak{D}$ and on  $[0,T]\times 
\partial\mathfrak{B}_0$, and 3)
the initial condition 
$$ \bm{\varphi}(0,\bar{\bm{x}})=\bar{\bm{x}} \quad\mbox{for}\quad \forall 
\bar{\bm{x}}\in \tilde{\mathfrak{D}}$$
holds.
\end{Definition}


As for the existence of the inverse mapping of $\bm{\varphi}(t, \cdot)$, we have

\begin{Proposition}\label{Prop.EulerDomain}
Let $\bm{\varphi} \in \bigcap_{\ell =1,2}C^{2-\ell}([0, T];C^{\ell}(\mathbf{Cl}\tilde{\mathfrak{D}}) $ be a vector field such that
$$ \bm{\varphi}(0,\bar{\bm{x}})=\bar{\bm{x}} \quad\mbox{for}\quad \forall 
\bar{\bm{x}}\in \tilde{\mathfrak{D}},$$
Let ${\mathfrak{O}}, \mathfrak{O}_0, \mathfrak{O}_1$ be connected open subsets of $\mathbb{R}^3$ such that
\begin{equation}
\mathfrak{B}_0 \subset\mathfrak{O}_0\Subset{\mathfrak{O}}
\Subset \mathfrak{O}_1
\Subset\tilde{\mathfrak{D}}
\end{equation}
and $\mathfrak{O}_1$ is convex. Then
for a sufficiently small ${T_1} (<T)$ there is a mapping $\bm{\psi}(t,\cdot): {\mathfrak{O}} \rightarrow \mathfrak{O}_1$, $t$ being $\in [0,{T_1}]$, such that
\begin{equation}
\bm{x}=\bm{\varphi}(t, \bm{\psi}(t,\bm{x}))
\end{equation}
for $\forall (t,\bm{x}) \in [0,{T_1}]\times {\mathfrak{O}}$. Moreover
\begin{equation}
\mathfrak{O}_0 \subset \bm{\psi}(t, {\mathfrak{O}})=
\{ \bm{\psi}(t,\bm{x} ) \  |\  \bm{x} \in {\mathfrak{O}} \}
\end{equation}
for $\forall t \in [0,{T_1}]$.
\end{Proposition}

Proof. Fixing $t$, we are going to solve the equation for unknown $\bar{\bm{x}}$
$$\bm{x}=\bm{\varphi}(t,\bar{\bm{x}})
$$
for gven $\bm{x}$. Writing 
$$F(\bar{\bm{x}})=\bm{x}+\bar{\bm{x}}-\bm{\varphi}(t,\bar{\bm{x}}),$$
we convert the equation to the fixed point problem
$$\bar{\bm{x}}=F(\bar{\bm{x}}). $$
Note 
\begin{align*}
&\|F(\bar{\bm{x}})-\bm{x}\|=\|\bar{\bm{x}}-\bm{\varphi}(t,\bar{\bm{x}})\|=
o(1),\\
& \| DF(\bar{\bm{x}})\|=\| I -D\bm{\varphi}(t,\bar{\bm{x}})\|=o(1)
\end{align*}
uniformly for $\bar{\bm{x}} \in \mathfrak{O}_1$ as $t \rightarrow 0$, since
$ \bm{\varphi}(0,\bar{\bm{x}})=\bar{\bm{x}}$. Therefore, if $\bm{x} \in 
{\mathfrak{O}}$, we see
$F(\bar{\bm{x}}) \in \mathfrak{O}_1$ for $\bar{\bm{x}} \in \mathfrak{O}_1$ and
$\|F(\bar{\bm{x}}')-F(\bar{\bm{x}})\|\leq \Theta \|\bar{\bm{x}}'-\bar{\bm{x}}\| $
for $\bar{\bm{x}}', \bar{\bm{x}} \in \mathfrak{O}_1$, $\Theta$ being $\in ]0,1[$, provided that $0\leq t \leq {T_1} \ll 1$. Here recall that $\mathfrak{O}_1$ is supposed to be convex. Therefore the fixed point problem admits the unique soluton $\bar{\bm{x}}=\bm{\psi}(t,\bm{x}) \in \mathfrak{O}_1$ for $(t,\bm{x})
\in [0,{T_1}]\times {\mathfrak{O}}$. Then
$$\bm{x}=\bm{\varphi}(t,\bm{\psi}(t,\bm{x}))
$$
for $(t,\bm{x})\in [0, {T_1}]\times{\mathfrak{O}}$.

Since $\mathbf{Cl}\mathfrak{O}_0 \subset {\mathfrak{O}}$, taking ${T_1}$ smaller if necessary, we can assume that
$\bm{\varphi}(t,\bar{\bm{x}}) \in \tilde{\mathfrak{D}}$ if
$(t, \bar{\bm{x}}) \in [0, {T_1}] \times
\mathfrak{O}_0$. Thus, if $\bar{\bm{x}}\in \mathfrak{O}_0$, there is $\bm{x} \in {\mathfrak{O}}$ such that
$\bar{\bm{x}}=\bm{\psi}(t,\bm{x})$ provided that $t \in [0,{T_1}]$. 
$\square$\\

In this situation, we can claim

\begin{Theorem}
Let $\mathfrak{D}$ be an admissible domain with cover $\tilde{\mathfrak{D}}$, $\overset{\circ}{\varUpsilon}$ be an admissible $\varUpsilon$-data
on $\mathfrak{D}$, and $\varphi$ be an admissible flow associated with $\overset{\circ}{\varUpsilon}$. Let 
${T_1}, {\mathfrak{O}}$ be those of Proposition \ref{Prop.EulerDomain}. 
Put
\begin{align}
&\varUpsilon(t,\bm{x})=
\overset{\circ}{\varUpsilon}(\bm{\psi}(t, \bm{x}))\Big(\mathrm{det}
D\bm{\varphi}(t, \bar{\bm{x}})\Big|_{\bar{\bm{x}}=\bm{\psi}(t,\bm{x})}\Big)^{-(\gamma-1)}  \\
&
\bm{v}(t,\bm{x})=\frac{\partial}{\partial t}\bm{\varphi}(t, \bar{\bm{x}})\Big|_{\bar{\bm{x}}=\bm{\psi}(t,\bm{x})}
\end{align}
for $(t, \bm{x}) \in [0,{T_1}]\times {\mathfrak{O}}$.
Then $(\varUpsilon,\bm{v})$ is a classical $(\varUpsilon,\bm{v})$-solution on $[0,{T_1}] \times
\underline{\mathfrak{D}}$, where $\underline{\mathfrak{D}}={\mathfrak{O}}\setminus
\mathfrak{B}_0$.
\end{Theorem}

\textbullet \  
Now let us derive the linearized approximation of \eqref{ECL},
\eqref{EM2}, \eqref{BCL}. Namely, we fix a stationary solution 
$(\bar{\Upsilon}, \bm{0})$,
say, we take 
\begin{equation}
\bar{\varUpsilon}(\bm{x})=-\Phi^{\langle\Omega\rangle}(\bm{x})
-\frac{\mathsf{G}M_0}{R} \label{3.24}
\end{equation}
under the assumption {\bf (K)}.
Let us denote 
\begin{equation}
\mathfrak{R}:=\{ \bar{\Upsilon} >0, \varpi < 3R/2, r>R_0 \}
=\{ \bar{\rho} >0 , r>R_0 \}.
\end{equation}

Considering small $\varepsilon$, a quantity $Q$ will denoted by $\bm{O}(\varepsilon)$ if $Q$ and its derivatives are of order $O(\varepsilon)$ uniformly on each bounded interval of $t$. We assume that $\Upsilon-\bar{\Upsilon}, \bm{v}$ are of $\bm{O}(\varepsilon)$. \\

Then \eqref{3.6} shows
\begin{align}
D_{\bar{\bm{x}}}\bm{x}(t,\bar{\bm{x}})&=I + \int_0^t(D_{\bm{x}}\bm{v}(s,\bm{\varphi}(s,\bar{\bm{x}}))ds+\bm{O}(\varepsilon^2) \\
&=I+ \bm{O}(\varepsilon), \\
D_{\bm{x}}\bar{\bm{x}}(t,\bm{x})&=I +\bm{O}(\varepsilon),
\end{align}
and so on. 

The equation \eqref{ECL} gives
\begin{equation}
\Upsilon^L(t,\bar{\bm{x}})=\overset{\circ}{\Upsilon}(\bar{\bm{x}})\Big(
1-(\gamma-1)\mathrm{div}_{\bar{\bm{x}}}(\bm{\varphi}(t,\bar{\bm{x}})-\bar{\bm{x}})\Big)+\bm{O}(\varepsilon^2).
\end{equation}

Recalling \eqref{3.24}, 
we have 
\begin{equation}
\Upsilon^L(t,\bar{\bm{x}})+\Phi^{\langle\Omega\rangle}(\bm{\varphi}(t,\bar{\bm{x}}))=
\Upsilon^L(t,\bar{\bm{x}})-\bar{\Upsilon}_L(\bar{\bm{x}})-\frac{\mathsf{G}M_0}{R}.
\end{equation}

Here we define $\bar{\Upsilon}_L$ by
\begin{equation}
\bar{\Upsilon}_L(\bar{\bm{x}})=\bar{\Upsilon}(\bm{\varphi}(t,\bar{\bm{x}})).
\end{equation}
Then
\begin{equation}
\bar{\Upsilon}_L(\bar{\bm{x}})=\bar{\Upsilon}(\bar{\bm{x}})+
(\mathrm{grad}_{\bar{\bm{x}}}\bar{\Upsilon}_L(\bar{\bm{x}})|\bm{\varphi}(t,\bar{\bm{x}})-\bar{\bm{x}}) +\bm{O}(\varepsilon^2).
\end{equation}

Summing up, putting
\begin{align}
G(t,\bar{\bm{x}})&:=\Upsilon(t,\bm{\varphi}(t,\bar{\bm{x}})+\Phi^{\langle\Omega\rangle}(\bm{\varphi}(t,\bar{\bm{x}})) + \frac{\mathsf{G}M_0}{R}, \nonumber \\
&=\Upsilon^L(t,\bar{\bm{x}})-\bar{\Upsilon}_L(\bar{\bm{x}}),
\end{align}
we see
\begin{align}
G(t,\bar{\bm{x}})&=\overset{\circ}{\Upsilon}(\bar{\bm{x}})-
\bar{\Upsilon}(\bar{\bm{x}})-(\gamma-1)\overset{\circ}{\Upsilon}(\bar{\bm{x}})
\mathrm{div}_{\bar{\bm{x}}}(\bm{\varphi}(t,\bar{\bm{x}})-\bar{\bm{x}}) + \nonumber \\
&-(\mathrm{grad}_{\bar{\bm{x}}}\bar{\Upsilon}_L(\bar{\bm{x}})|\bm{\varphi}(t,\bar{\bm{x}})-\bar{\bm{x}})+
\bm{O}(\varepsilon^2).
\end{align}\\

Let a constant vector field $\overset{\circ}{\bm{\xi}}(\bar{\bm{x}})=\bm{O}(\varepsilon)$ on $\mathfrak{R}$ be given arbitrarily. We take $\overset{\circ}{\varUpsilon}$ given by
\begin{equation}
\overset{\circ}{\varUpsilon}(\bar{\bm{x}})=\bar{\varUpsilon}(\bar{\bm{x}})
-(\gamma-1)\bar{\varUpsilon}(\bar{\bm{x}})\mathrm{div}_{\bar{\bm{x}}}\overset{\circ}{\bm{\xi}}
-(\mathrm{grad}_{\bar{\bm{x}}}\bar{\varUpsilon}_L(\bar{\bm{x}})|\overset{\circ}{\bm{\xi}}).
\label{IUp}
\end{equation}

Then it turns out to be that
\begin{align}
G(t,\bar{\bm{x}})=-(\gamma-1)\bar{\Upsilon}(\bar{\bm{x}})\mathrm{div}_{\bar{\bm{x}}}\hat{\bm{\xi}}-
(\mathrm{grad}_{\bar{\bm{x}}}\bar{\Upsilon}_L(\bar{\bm{x}})\Big| \hat{\bm{\xi}})+\bm{O}(\varepsilon^2),\\
=-(\gamma-1)\bar{\Upsilon}(\bar{\bm{x}})\mathrm{div}_{\bar{\bm{x}}}\hat{\bm{\xi}}-
(\mathrm{grad}_{\bar{\bm{x}}}\bar{\Upsilon}(\bar{\bm{x}})\Big| \hat{\bm{\xi}})+\bm{O}(\varepsilon^2),
\end{align}
where
\begin{equation}
\hat{\bm{\xi}}(t,\bar{\bm{x}})=\bm{\varphi}(t,\bar{\bm{x}})-\bar{\bm{x}}+\overset{\circ}{\bm{\xi}}(\bar{\bm{x}})
\end{equation}
and
$\mathrm{grad}_{\bar{\bm{x}}}\bar{\varUpsilon}(\bar{\bm{x}})$ means
$\Big(\mathrm{grad}_{\bm{x}}\bar{\varUpsilon}(\bm{x})\Big)\Big|_{\bm{x}=\bar{\bm{x}}}$ so that
\begin{align*}
\mathrm{grad}_{\bar{\bm{x}}}\bar{\varUpsilon}_L(\bar{\bm{x}})&=(D_{\bar{\bm{x}}}\bm{x}(t,\bar{\bm{x}}))^{\top}\mathrm{grad}_{\bar{\bm{x}}}\bar{\varUpsilon}(\bar{\bm{x}}) \\
&=\mathrm{grad}_{\bar{\bm{x}}}\bar{\varUpsilon}(\bar{\bm{x}})+\bm{O}(\varepsilon).
\end{align*}

Then we see
\begin{equation}
G(t,\bar{\bm{x}})
=-\bar{\sigma}(\bar{\bm{x}})
\mathrm{div}_{\bar{\bm{x}}}(\bar{\rho}(\bar{\bm{x}})\hat{\bm{\xi}}(t,\bar{\bm{x}}))+
\bm{O}(\varepsilon^2) \quad\mbox{on}\quad \mathfrak{R},
\end{equation}
where
\begin{equation}
\sigma=\frac{d\Upsilon}{d\rho}=\mathsf{A}\gamma \rho^{\gamma-2}.
\end{equation}

Moreover we introduce
\begin{equation}
\bm{\xi}(t,\bar{\bm{x}}):=
(D\bm{\varphi}(t,\bar{\bm{x}}))^{-1} 
\hat{\bm{\xi}}(t,\bar{\bm{x}}) 
-\int_0^t
\frac{\partial }{\partial  s}\Big[ (D\bm{\varphi}(s,\bar{\bm{x}}))^{-1}\Big]
\hat{\bm{\xi}}(s,\bar{\bm{x}})ds.
\end{equation}
Then we have
\begin{equation}
\bm{\xi} =\hat{\bm{\xi}}+\bm{O}(\varepsilon^2),
\end{equation}

\begin{equation}
\frac{\partial  \bm{\xi} }{\partial  t}=\bm{v}^L, \label{LgV}
\end{equation}
\begin{equation}
\bm{\xi}(0,\bar{\bm{x}})=\overset{\circ}{\bm{\xi}}(\bar{\bm{x}})
\end{equation}
and
\begin{equation}
G(t,\bar{\bm{x}})=
-\bar{\sigma}(\bar{\bm{x}})
\mathrm{div}_{\bar{\bm{x}}}(\bar{\rho}(\bar{\bm{x}})\bm{\xi}(t,\bar{\bm{x}}))+\bm{O}(\varepsilon^2) \label{FinG}
\end{equation}\\

The approximation of the equation of motion \eqref{EM2} is clearly given by
\begin{equation}
\frac{\partial  \bm{v}^L}{\partial  t}+B\bm{v}^L+\mathrm{grad}_{\bar{\bm{x}}}G=\bm{O}(\varepsilon^2).
\label{LinEq2}
\end{equation}

Inserting \eqref{LgV}, \eqref{FinG} into \eqref{LinEq2}, we have
\begin{equation}
\frac{\partial ^2\bm{\xi}}{\partial  t^2}+
B\frac{\partial  \bm{\xi}}{\partial  t}+
\mathrm{grad}_{\bar{\bm{x}}}\Big(
-\bar{\sigma}(\bar{\bm{x}})
\mathrm{div}_{\bar{\bm{x}}}(\bar{\rho}(\bar{\bm{x}})\bm{\xi})\Big)=\bm{O}(\varepsilon^2).
\end{equation}\\

The approximation of the boundary condition \eqref{BCL}
is given by
\begin{equation}
\Big(\bm{n}(\bar{\bm{x}})\Big|\frac{\partial \bm{\xi}}{\partial t}
(t,\bar{\bm{x}})\Big)=0
\quad\mbox{for}\quad (t,\bar{\bm{x}})\in [0,T]\times
\partial \mathfrak{B}_0. \label{approxBC1}
\end{equation}

In fact we have
$$\bm{n}^L=
\frac{1}{R_0}(I+\bm{O}(\varepsilon))(\bar{\bm{x}}+\bm{O}(\varepsilon))=
\bm{n}(\bar{\bm{x}})+\bm{O}(\varepsilon)$$
and 
$$\bm{v}^L=\frac{\partial \bm{\xi}}{\partial t}=\bm{O}(\varepsilon).$$
Of course \eqref{approxBC1} is equivalent  to
\begin{equation}
(\bm{n}(\bar{\bm{x}})| \bm{\xi}(t,\bar{\bm{x}}) )=0
\quad\mbox{for}\quad (t,\bar{\bm{x}})\in [0,T]\times
\partial \mathfrak{B}_0, \label{approxBC}
\end{equation}
provided that
\begin{equation}
(\bm{n}(\bar{\bm{x}})|\overset{\circ}{\bm{\xi}} (\bar{\bm{x}}) )=0 \quad\mbox{on}\quad
\partial \mathfrak{B}_0.
\end{equation}

====\\

Now the domain $\mathfrak{R}=
\{\bar{\rho}>0\}=\{\bar{\varUpsilon}>0\}\cap \{\varpi <3R/2\}$
has the form
$$\mathfrak{R}=\{ R_0<r <R\cdot H(\zeta^2;\kappa,\lambda),\quad  -1 \leq \zeta:=\frac{z}{r}\leq 1 \},$$
where $\zeta^2 \mapsto H(\zeta^2;\kappa, \lambda)$ is a smooth monotone function on $[0,1]$ such that
$1=H(1;\kappa,\lambda)\leq H(0;\kappa, \lambda)$. (Note $H(1)=H(0) \Leftrightarrow \kappa=0, \mbox{that is,} \Omega=0$.)
Given a small $\overset{\circ}{\bm{\xi}}$, we consider 
$\overset{\circ}{\varUpsilon}$ determined by \eqref{IUp}.
Of course $\overset{\circ}{\varUpsilon}-\bar{\varUpsilon}$ is small, but
the topology of 
$\{ \overset{\circ}{\varUpsilon} > 0\}\cap \{\varpi <3R/2\}$
is not clear, generaly speaking. In fact, as for
$ \{ \overset{\circ}{\varUpsilon} >0 \} = 
\Big\{\Big(1-(\gamma-1)(\nabla|\overset{\circ}{\bm{\xi}})\Big)\bar{\varUpsilon} >
(\nabla \bar{\varUpsilon}|\overset{\circ}{\bm{\xi}}) \Big\}$,
we have that $\displaystyle \frac{\nabla\bar{\varUpsilon}}{\bar{\varUpsilon}}$ may diverge along the vacuum boundary
 $\Sigma_1=\{\bar{\varUpsilon}=0, \varpi <\frac{3R}{2}\}=\{\bar{\rho}=0 \}
$ of $\bar{\rho}$. At least, if, for example, $\overset{\circ}{\mbox{\boldmath$\xi$}}_1 \in C_0^{\infty}(\Sigma_0\cup \mathfrak{R})$, then for  $\overset{\circ}{\mbox{\boldmath$\xi$}}=
\varepsilon\overset{\circ}{\mbox{\boldmath$\xi$}}_1$ with $\varepsilon \ll 1$, it is guaranteed that $\{ \overset{\circ}{\varUpsilon} >0\}\cap\{ \varpi <3R/2\}=\mathfrak{R}
(=\{ \bar{\varUpsilon} > 0 \}\cap\{ \varpi <3R/2\})
$.
More generally, {\bf if
$\displaystyle
\frac{1}{\bar{\rho}}(\nabla|\bar{\rho}\overset{\circ}{\bm{\xi}}_1)=
 (\nabla\log\bar{\rho}|\overset{\circ}{\bm{\xi}}_1)+(\nabla|\overset{\circ}{\bm{\xi}}_1) \subset
 \frac{1}{\gamma-1}(\nabla\log\bar{\varUpsilon}|\overset{\circ}{\bm{\xi}}_1)+(\nabla|\overset{\circ}{\bm{\xi}}_1)$
is bounded on a neighborhood of the vacuum boundary $\Sigma_1$ of $\bar{\rho}$, then $\bar{\rho}>0 \Leftrightarrow
\overset{\circ}{\varUpsilon}>0$ there for $\overset{\circ}{\bm{\xi}}=\varepsilon\overset{\circ}{\bm{\xi}}_1$ with $\varepsilon \ll 1$}.\\

{\bf Historical Remark:}\   The derivation of the linearized approximation of the equations in Lagrangian co-ordinate system can be found
 \cite[Sect. 56]{LedouxW}, \cite[pp. 139-140.]{Batchelor}, \cite[p.11, (A)]{Bjerknes}, \cite{LyndenBO}, \cite[p.500, (1)]{Lebovitz} and so on. But there was considered only the case of $\overset{\circ}{\rho}=\bar{\rho}$
 and $\overset{\circ}{\bm{\xi}}=0$.\\

\section{Linearized equations for perturbations from a static stationary solution}

Let us fix a static stationary solution
\begin{equation}
\bar{\varUpsilon}=
-\Phi^{\langle \Omega \rangle}-\frac{{\mathsf{G}M_0}}{R},
\quad \bar{\bm{v}}=\mathbf{0}
\end{equation}
under the assumption {\bf (K)}. We are concerned with the domain
$$\mathfrak{R}=\{ \bar{\rho} >0\}=\{ \bar{\varUpsilon} >0, \varpi <\frac{3R}{2}, R_0 <r \}. $$
We consider the initial boundary problem
\begin{subequations}
\begin{align}
&\frac{\partial^2 \bm{\xi}}{\partial t^2}+
\bm{B}\frac{\partial \bm{\xi} }{\partial t}+
\bm{L}\bm{\xi}=0 &\quad\mbox{on}\quad [0,+\infty[\times \mathfrak{R}, \label{LEAMa} \\
&(\bm{n}|\bm{\xi})=0 &\quad\mbox{on}\quad [0,+\infty[\times \Sigma_0, \label{LEAMb}\\
&\bm{\xi}\Big|_{t=0}=\overset{\circ}{\bm{\xi}}(\bm{x}),\quad 
\frac{\partial\bm{\xi}}{\partial t}\Big|_{t=0}=\overset{\circ}{\bm{v}}(\bm{x}) &\quad\mbox{on}\quad \mathfrak{R}. \label{LEAMc}
\end{align}
\end{subequations}
Here 
\begin{align}
&\bm{B}\bm{v}=2\Omega\bm{J}\bm{v},\quad \bm{J}=
\begin{bmatrix}
0 & -1 & 0 \\
1 & 0 & 0 \\
0 & 0 & 0
\end{bmatrix}, \\
&\mbox{\boldmath$L$}\mbox{\boldmath$\xi$}=\nabla G, \quad G=-\overline{\frac{d\varUpsilon}{d\rho}}g,\quad g=(\nabla | \bar{\rho}\mbox{\boldmath$\xi$}), \\
&\bm{n}=\frac{\bm{x}}{R_0}
\end{align}
and $\overset{\circ}{\bm{\xi}}, \overset{\circ}{\bm{v}}$ are given vector fields on $\mathfrak{R}$ such that $(\bm{n}|\overset{\circ}{\bm{\xi}})=(\bm{n}|\overset{\circ}{\bm{v}})=0$ on $\Sigma_0$, while
$$\Sigma_0=\partial\mathfrak{B}_0=\{ r= R_0\}.$$

Here and hereafter the Lagrangian coordinate, denoted by $\bar{\bm{x}}$ in thr previous section, is denoted by $\bm{x}$, since we do not refer the Eulerian coordinate so for.\\

We consider the differential operator $\mbox{\boldmath$L$}$ in the Hilbert space 
$\mathfrak{H}$ of all measurable functions $\mbox{\boldmath$\xi$}$
defined on $\mathfrak{R}$ such that
$\|\mbox{\boldmath$\xi$}\|_{\mathfrak{H}} < \infty$, where
\begin{equation}
\|\mbox{\boldmath$\xi$}\|_{\mathfrak{H}}^2=\int_{\mathfrak{R}}\|\mbox{\boldmath$\xi$}(\mbox{\boldmath$x$})\|^2\bar{\rho}(\mbox{\boldmath$x$})d\mbox{\boldmath$x$},
\end{equation}
that is, the inner product of $\mathfrak{H}$ is
\begin{equation}
(\mbox{\boldmath$\xi$}_1|\mbox{\boldmath$\xi$}_2)_{\mathfrak{H}}=
\int_{\mathfrak{R}}(\mbox{\boldmath$\xi$}_1(\mbox{\boldmath$x$})|\mbox{\boldmath$\xi$}_2(\mbox{\boldmath$x$}))\bar{\rho}(\mbox{\boldmath$x$})d\mbox{\boldmath$x$}.
\end{equation}
Of course
$$(\mbox{\boldmath$\xi$}_1(\mbox{\boldmath$x$})|\mbox{\boldmath$\xi$}_2(\mbox{\boldmath$x$})):=\sum_k \xi_1^k(\mbox{\boldmath$x$})(\xi_2^k(\mbox{\boldmath$x$}))^*
\quad\mbox{for}\quad 
\mbox{\boldmath$\xi$}_{\mu}(\mbox{\boldmath$x$})=
\begin{bmatrix}
\xi_{\mu}^1(\mbox{\boldmath$x$}) \\
\\
\xi_{\mu}^2(\mbox{\boldmath$x$})\\
\\
 \xi_{\mu}^ 3(\mbox{\boldmath$x$})
\end{bmatrix}, \mu=1,2.$$\\

 Here and hereafter $Z^*$ denotes the complex conjugate $X-\mathrm{i}Y$of $Z=X+\mathrm{i}Y$, while $\mathrm{i}$ stands for the imaginary unit, $\sqrt{-1}$.\\

Briefly speaking, we consider $\mathfrak{H}=L^2(\mathfrak{R}, \bar{\rho}dx; \mathbb{C}^3)$.\\

For $\mbox{\boldmath$\xi$}_{\mu}\in C_0^{\infty}(\mathfrak{R}), \mu=1,2$, we have
\begin{equation}
(\mbox{\boldmath$L$}\mbox{\boldmath$\xi$}_1|\mbox{\boldmath$\xi$}_2)_{\mathfrak{H}}=
\int_{\mathfrak{R}}
\overline{\frac{d\varUpsilon}{d\rho}}g_1g_2^*d\mbox{\boldmath$x$},\quad\mbox{with}\quad
g_{\mu}=(\nabla | \bar{\rho}\mbox{\boldmath$\xi$}_{\mu}).
\end{equation}\\

So we put
\begin{equation}
Q(\mbox{\boldmath$\xi$}_1,\mbox{\boldmath$\xi$}_2)=
\int_{\mathfrak{R}}
\overline{\frac{d\varUpsilon}{d\rho}}g_1g_2^*d\mbox{\boldmath$x$}\quad\mbox{with}\quad
g_{\mu}=(\nabla| \bar{\rho}\mbox{\boldmath$\xi$}_{\mu}),
\end{equation}
and
\begin{equation}
Q[\mbox{\boldmath$\xi$}]=Q(\mbox{\boldmath$\xi$}, \mbox{\boldmath$\xi$})=\int_{\mathfrak{R}}\overline{\frac{d\varUpsilon}{d\rho}}|g|^2d\mbox{\boldmath$x$}\quad 
\mbox{with}\quad g=(\nabla | \bar{\rho}\mbox{\boldmath$\xi$}). \label{3.19}
\end{equation}\\

We start from the operator $\mathfrak{T_1}_c$ in $\mathfrak{H}$ defined by
 $\mathfrak{T_1}_c : \mbox{\boldmath$\xi$} \mapsto \mbox{\boldmath$L$}\mbox{\boldmath$\xi$} $ on the domain 
$\mathsf{D}(\mathfrak{T_1}_c)=C_0^{\infty}(\mathfrak{R})$. Then $\mathfrak{T_1}_c$ is
densely defined,  symmetric and bounded from below as
$$(\mathfrak{T_1}_c\mbox{\boldmath$\xi$}|\mbox{\boldmath$\xi$})_{\mathfrak{H}}=Q[\mbox{\boldmath$\xi$}] \geq 0.
$$Therefore, seeing, e.g., \cite[Chapter VI, Section 2.3]{Kato}, we have that $\mathfrak{T_1}_c$ admits the Friedrichs extension $\mathfrak{T_1}$ which is a self-adjoint operator in $\mathfrak{H}$.  The domain of $\mathfrak{T_1}$ is 
\begin{equation}
\mathsf{D}(\mathfrak{T_1})=\{ \mbox{\boldmath$\xi$} \in \mathfrak{H}_0^{\mathrm{div}}\  |\  \mbox{\boldmath$L$}\mbox{\boldmath$\xi$} \in \mathfrak{H} \quad \mbox{in the sense of distribution} \}.
\end{equation}
Here $\mathfrak{H}_0^{\mathrm{div}}$ is the Hilbert space of all 
$\mbox{\boldmath$\xi$} \in \mathfrak{H}$ such that there is a sequence
$\mbox{\boldmath$\varphi$}_n \in C_0^{\infty}(\mathfrak{R})$ such that
\begin{align*}
&\mbox{\boldmath$\varphi$}_n \rightarrow \mbox{\boldmath$\xi$} \quad\mbox{in}\quad \mathfrak{H}\quad\mbox{as}\quad
n \rightarrow \infty, \\
&Q[\mbox{\boldmath$\varphi$}_m-\mbox{\boldmath$\varphi$}_n] \rightarrow 0\quad\mbox{as}\quad
n,m \rightarrow \infty.
\end{align*}

In order to fix the idea we use
\begin{Definition}
We put
\begin{equation}
\mathfrak{H}^{\mathrm{div}}=\{ \bm{\xi} \in \mathfrak{H} \  |\  
g=(\nabla|\bar{\rho}\bm{\xi}) \in L^2(\mathfrak{R}; 
\overline{\frac{d\varUpsilon}{d\rho}}d\bm{x})
 \}, \label{DefHdiv}
\end{equation}
and regard it as a Hilbert space endowed with the inner product
\begin{equation}
(\mbox{\boldmath$\xi$}_1|\mbox{\boldmath$\xi$}_2)_{\mathfrak{H}^{\mathrm{div}}}=
Q(\mbox{\boldmath$\xi$}_1, \mbox{\boldmath$\xi$}_2)+(\mbox{\boldmath$\xi$}_1|\mbox{\boldmath$\xi$}_2)_{\mathfrak{H}}.
\end{equation}
\end{Definition}

Thus we are saying

\begin{Definition}  $\mathfrak{H}_0^{\mathrm{div}}$ is the closure of $C_0^{\infty}(\mathfrak{R})$ in $\mathfrak{H}^{\mathrm{div}}$.
\end{Definition}

Hereafter we shall denote by $\mbox{\boldmath$L$}$ the Friedrichs extension $\mathfrak{T_1}$, diverting the letter. So we have

\begin{Theorem}
The operator $\mbox{\boldmath$L$}$ is a self-adjoint operator bounded from below by $0$ in the Hilbert
space $\mathfrak{H}$, whose domain is 
\begin{equation}
\mathsf{D}(\mbox{\boldmath$L$})=\{ \mbox{\boldmath$\xi$} \in \mathfrak{H}_0^{\mathrm{div}} \  |\   \mbox{\boldmath$L$}\mbox{\boldmath$\xi$} \in \mathfrak{H}\}.
\end{equation}
\end{Theorem}

Note that $\mbox{\boldmath$\xi$} \in \mathsf{D}(\mbox{\boldmath$L$}) \subset \mathfrak{H}_0^{\mathrm{div}}$ enjoys the boundary condition
\begin{equation}
(\bm{n}|\bm{\xi}) =0
\quad\mbox{on}\quad
\Sigma_0=\{ r=R_0\} \label{0416}
\end{equation}
in the following sense: There is known to uniquely exist the `normal trace operator' 
$\gamma_{\mathrm{n}}$ which maps $\mathfrak{H}^{\mathrm{div}}$ into $H^{-1/2}(\Sigma_0)$ continuously such that $\gamma_{\mathrm{n}}\mbox{\boldmath$\xi$}=-(\bm{n}|\mbox{\boldmath$\xi$}) \in C^1(\mathfrak{R}\cup\partial\mathfrak{R})$, and $\gamma_{\mathrm{n}}(\mbox{\boldmath$\xi$})=0$ when 
$\mbox{\boldmath$\xi$} \in \mathfrak{H}_0^{\mathrm{div}}$. (See e.g., \cite[Chapter I]{GR86}.)
Here we have used the fact that $\displaystyle \bar{\rho}, \overline{\frac{d\varUpsilon}{d\rho}} \in C^{\infty}(\mathfrak{R}\cup\Sigma_0)$ 
and
\begin{align*}
&\bar{\rho}\geq
\Big(\frac{\gamma-1}{\mathsf{A}\gamma}\Big)^\frac{1}{\gamma-1}
\Big({\mathsf{G}M_0}\Big(\frac{1}{R_0}-\frac{1}{R}\Big)\Big)^{\frac{1}{\gamma-1}}>0
\quad\mbox{on}\  \Sigma_0, \\
&\overline{\frac{d\varUpsilon}{d\rho}}
 \geq
(\gamma-1)\Big(\frac{\mathsf{A}\gamma}{\gamma-1}\Big)^{\frac{1}{\gamma-1}}
\Big({\mathsf{G}M_0}\Big(\frac{1}{R_0}-\frac{1}{R}\Big)\Big)^{-\frac{2-\gamma}{\gamma-1}}
>0 
\quad\mbox{on}\  \Sigma_0,
\end{align*}
thanks to the assumption $R>R_0$,
so that 
\begin{align*}
&\|\mbox{\boldmath$\xi$}\|^2\bar{\rho}d\mbox{\boldmath$x$} =
\frac{1}{\bar{\rho}}\|\hat{\mbox{\boldmath$\xi$}}\|^2d\mbox{\boldmath$x$}
\geq \frac{1}{C}
\|\hat{\mbox{\boldmath$\xi$}}\|^2d\mbox{\boldmath$x$}\quad  \mbox{on}\  \mathfrak{R}, \\
&\overline{\frac{d\varUpsilon}{d\rho}}|(\nabla|\bar{\rho}\mbox{\boldmath$\xi$})|^2d\mbox{\boldmath$x$} =
{\mathsf{A}\gamma}{\bar{\rho}^{\gamma-2}}
|(\nabla|\hat{\mbox{\boldmath$\xi$}})|^2d\mbox{\boldmath$x$}
\geq
\frac{1}{C}
|(\nabla|\hat{\mbox{\boldmath$\xi$}})|^2d\mbox{\boldmath$x$}
\quad  \mbox{on}\  \mathfrak{R},
\end{align*}
where $\hat{\mbox{\boldmath$\xi$}}:=
\bar{\rho}{\mbox{\boldmath$\xi$}}$
and $C$ is a sufficiently large finite positive number, and the theory on the functional spaces $H({\mathrm{div}}; \mathfrak{R}), H_0({\mathrm{div}};\mathfrak{R}), H^{1/2}(\partial\mathfrak{R}), H^{-1/2}(\partial\mathfrak{R})$ and the normal trace operator $\gamma_{\mathrm{n}}$ seen in \cite{GR86} can be applied
to the vector filed $\hat{\mbox{\boldmath$\xi$}}=\bar{\rho}\mbox{\boldmath$\xi$}$,
since the above estimates guarantee that the space 
$\hat{\mathfrak{H}}^{\mathrm{div}}:=\{ \hat{\mbox{\boldmath$\xi$}} | \mbox{\boldmath$\xi$} \in \mathfrak{H}^{\mathrm{div}}\},
\hat{\mathfrak{H}}_0^{\mathrm{div}}:=\{ \hat{\mbox{\boldmath$\xi$}} | \mbox{\boldmath$\xi$} \in \mathfrak{H}_0^{\mathrm{div}}\}
$ are continuously imbedded into 
$H(\mathrm{div}; \mathfrak{R}),
H_0(\mathrm{div}; \mathfrak{R})$. 

Therefore, if $\bm{\xi}(t,\cdot) \in \mathsf{D}(\bm{L})\quad \forall t$, then
the boundary condition \eqref{LEAMb} is satisfied in this sense.\\ 

As for the operator $\mbox{\boldmath$J$}\cdot : \bm{v} \mapsto \mbox{\boldmath$J$}\bm{v}$, it is clear that it is a bounded linear operator on $\mathfrak{H}$. But, for the sake of generality, we replace $2\Omega\bm{J}$ by the operator $\bm{B}:\bm{v}\mapsto B\bm{v}$ defined by an arbitrary smooth matrix-valued function $B : \bm{x} \mapsto
(B_j^i(\bm{x}))_{i,j} \in  C^{\infty}(\mathfrak{R}\cup \partial\mathfrak{R} ;  \mathbb{R}^{3\times 3})$. Of course $\bm{B}$ is a bounded operator on $\mathfrak{H}$, and
\begin{equation}
|(\bm{B}\bm{v}|\bm{v})_{\mathfrak{H}}| \leq \beta \|\bm{v}\|_{\mathfrak{H}},
\end{equation}
where $\beta=\|B(\cdot)\|_{L^{\infty}}$. The adjoint operator $\bm{B}^*$ is defined by the transposed matrix $B(\bm{x})^{\top}$. 
Later we shall use the following densely defined closed operator
\begin{equation}
\bm{L}^{c\bm{B}}:=\bm{L}+c\bm{B},\quad \mathsf{D}(\bm{L}^{c\bm{B}})=\mathsf{D}(\bm{L}),
\end{equation}
$c$ being a real number. When $c\not=0$, we cannot say that it is self-adjoint, in genaral, but we can claim the folowing

\begin{Proposition}\label{Resolv}
For any $ c \in \mathbb{R}$ and $\lambda >|c\beta|$ the operator
 $\mbox{\boldmath$L$} +c\bm{B}+\lambda$ has the bounded linear inverse operator $(\mbox{\boldmath$L$}+c\bm{B}+\lambda)^{-1}$ defined on the whole space $\mathfrak{H}$
such that
\begin{equation}
|\|(\mbox{\boldmath$L$}+c\bm{B}+\lambda)^{-1}\||_{\mathcal{B}(\mathfrak{H})}\leq \frac{1}{\lambda-|c\beta|}.
\end{equation}
\end{Proposition}

Proof. First we see $\mbox{\boldmath$L$} +c\bm{B}+\lambda$ is invertible. In fact, if
$$(\mbox{\boldmath$L$} +c\bm{B}+\lambda)\mbox{\boldmath$\xi$}=\mbox{\boldmath$f$},\quad \mbox{\boldmath$\xi$} \in \mathsf{D}(\mbox{\boldmath$L$}),\quad \mbox{\boldmath$f$} \in \mathfrak{H},
$$
then  we see 
$$
(\lambda-|c\beta|)\|\mbox{\boldmath$\xi$}\|_{\mathfrak{H}}^2 \leq Q[\mbox{\boldmath$\xi$}]
+(c\bm{B}\mbox{\boldmath$\xi$}|\mbox{\boldmath$\xi$})_{\mathfrak{H}}
+\lambda\|\mbox{\boldmath$\xi$}\|_{\mathfrak{H}}^2 = (\mbox{\boldmath$\xi$}|\mbox{\boldmath$f$})_{\mathfrak{H}} \leq \|\mbox{\boldmath$\xi$}\|_{\mathfrak{H}}\|\mbox{\boldmath$f$}\|_{\mathfrak{H}},
$$
therefore we have
$$\|\mbox{\boldmath$\xi$}\|_{\mathfrak{H}}\leq \frac{1}{\lambda-|c\beta|}\|\mbox{\boldmath$f$}\|_{\mathfrak{H}}.
$$
We claim that the range $\mathsf{R}(\mbox{\boldmath$L$}+c\bm{B}
+\lambda)$ is dense in $\mathfrak{H}$. In fact, suppose
$$((\mbox{\boldmath$L$}+c\bm{B}
 +\lambda)\mbox{\boldmath$\xi$}|\mbox{\boldmath$f$})=0\quad \forall \mbox{\boldmath$\xi$} \in \mathsf{D}(\mbox{\boldmath$L$}).
$$
Then 
\begin{align*}
(\mbox{\boldmath$L$}\mbox{\boldmath$\xi$}|\mbox{\boldmath$f$})&=
-((c\bm{B}+\lambda)\bm{\xi}|\bm{f}) \\
&=(\bm{\xi}|(-c\bm{B}^*-\lambda)\bm{f})
\end{align*}
for $\forall \bm{\xi} \in \mathsf{D}(\bm{L})$. 
Hence $\mbox{\boldmath$f$} \in \mathsf{D}(\mbox{\boldmath$L$}^*)$ and
$$\mbox{\boldmath$L$}^*\mbox{\boldmath$f$}=
-c\bm{B}^*-\lambda\mbox{\boldmath$f$}.$$
Since $\mbox{\boldmath$L$}=\mbox{\boldmath$L$}^*$, this means that $\mbox{\boldmath$f$} \in \mathsf{D}(\mbox{\boldmath$L$})$ and
$$(\mbox{\boldmath$L$}+c\bm{B}^*+\lambda)\mbox{\boldmath$f$}=0.$$
Since $\mbox{\boldmath$L$}+c\bm{B}^*+\lambda$ is invertible, for
$\|B^{\top}\|=\|B\|=\beta$, it follows that $\mbox{\boldmath$f$}=0$. Summing up, we have the assertion. $\square$\\

We are going to apply the Hille-Yosida theory to the initial-boundary value problem
\begin{align}
&\frac{\partial^2\mbox{\boldmath$\xi$}}{\partial t^2}+
\bm{B}
\frac{\partial\mbox{\boldmath$\xi$}}{\partial t}+
\mbox{\boldmath$L$}\mbox{\boldmath$\xi$}=\bf{0}, \nonumber \\
&\mbox{\boldmath$\xi$}=\overset{\circ}{\bm{\xi}},\quad \frac{\partial\mbox{\boldmath$\xi$}}{\partial t}=\overset{\circ}{\mbox{\boldmath$v$}}=
\begin{bmatrix}
v^1(0,\mbox{\boldmath$x$}) \\
\\
v^2(0,\mbox{\boldmath$x$}) \\
\\
v^3(0,\mbox{\boldmath$x$})
\end{bmatrix}
\quad \mbox{at}\quad t=0, \nonumber \\
&\mbox{\boldmath$\xi$}(t,\cdot)\in \mathsf{D}(\mbox{\boldmath$L$})\quad\mbox{for}\quad \forall t \geq 0. \label{IBV}
\end{align}

We put
\begin{align}
&U=
\begin{bmatrix}
\mbox{\boldmath$\xi$} \\
\\
\dot{\mbox{\boldmath$\xi$}}
\end{bmatrix},
\quad \dot{\mbox{\boldmath$\xi$}}=\frac{\partial\mbox{\boldmath$\xi$}}{\partial t}, \\
&\mathbf{A}U=
\begin{bmatrix}
O & -I \\
\\
\mbox{\boldmath$L$} & \bm{B}
\end{bmatrix}
U
=
\begin{bmatrix}
-\dot{\mbox{\boldmath$\xi$}} \\
\\
\bm{B}\dot{\mbox{\boldmath$\xi$}}+\mbox{\boldmath$L$}\mbox{\boldmath$\xi$}
\end{bmatrix}, \\
&\mathfrak{E}=\mathfrak{H}_0^{\mathrm{div}}\times \mathfrak{H} \nonumber \\
\mbox{with}& \nonumber \\
&(U_1|U_2)_{\mathfrak{E}}=
(\mbox{\boldmath$\xi$}_1|\mbox{\boldmath$\xi$}_2)_{\mathfrak{H}^{\mathrm{div}}}+(\dot{\mbox{\boldmath$\xi$}}_1|\dot{\mbox{\boldmath$\xi$}}_2)_{\mathfrak{H}}= \nonumber \\
&=Q(\mbox{\boldmath$\xi$}_1, \mbox{\boldmath$\xi$}_2)+(\mbox{\boldmath$\xi$}_1|\mbox{\boldmath$\xi$}_2)_{\mathfrak{H}}+
(\dot{\mbox{\boldmath$\xi$}}_1|\dot{\mbox{\boldmath$\xi$}}_2)_{\mathfrak{H}}, \\
&\mathsf{D}(\mathbf{A})=\mathsf{D}(\mbox{\boldmath$L$})\times \mathfrak{H}_0^{\mathrm{div}}.
\end{align}

Then the initial-boundary value problem \eqref{IBV} can be written as 
\begin{equation}
\frac{dU}{dt}+\mathbf{A}U=\mathbf{0},\quad U|_{t=0}=U_0, \label{IV}
\end{equation}
where 
\begin{equation}
U_0=
\begin{bmatrix}
\overset{\circ}{\bm{\xi}} \\
\\
\overset{\circ}{\mbox{\boldmath$v$}}
\end{bmatrix}.
\end{equation}

Applying \cite[Theorem 7.4]{Brezis}, we can claim
\begin{Proposition}\label{Prop.2}
If $U_0 \in \mathsf{D}(\mathbf{A})$, say, if
$\overset{\circ}{\bm{\xi}} \in \mathsf{D}(\bf{L})$ and 
$\overset{\circ}{\mbox{\boldmath$v$}} \in \mathfrak{H}_0^{\mathrm{div}}$, then there exists a unique solution 
$U \in C^1([0,+\infty[, \mathfrak{E})\cap C([0,+\infty[, \mathsf{D}(\mathbf{A}))
$
to the problem \eqref{IV}. Moreover $E(t)=\|U(t)\|_{\mathfrak{E}}^2$ enjoys
\begin{equation}
E(t)\leq e^{2\Lambda t}E(0),
\end{equation}
where $\Lambda=1 \vee \beta$.
\end{Proposition}

Here we consider that $\mathsf{D}(\mathbf{A})$ is equipped with the operator norm
$(\|U\|_{\mathfrak{E}}^2+\|\mathbf{A}U\|_{\mathfrak{E}}^2)^{1/2}$.

Proof of Proposition \ref{Prop.2}. Firstly $\mathbf{A}+1\vee \beta$ is monotone, that is, for $\forall U \in \mathsf{D}(\mathbf{A})$ we have
\begin{align*}
\mathfrak{Re}[(\mathbf{A}U|U)_{\mathfrak{E}}]+(1\vee \beta)\|U\|_{\mathfrak{E}}^2&=\mathfrak{Re}\Big[-Q(\dot{\mbox{\boldmath$\xi$}},\mbox{\boldmath$\xi$})-(\dot{\mbox{\boldmath$\xi$}}|\mbox{\boldmath$\xi$})_{\mathfrak{H}}+(\mbox{\boldmath$L$}\mbox{\boldmath$\xi$}|\dot{\mbox{\boldmath$\xi$}})_{\mathfrak{H}}+
(\mbox{\boldmath$B$}\dot{\mbox{\boldmath$\xi$}}|\dot{\mbox{\boldmath$\xi$}})_{\mathfrak{H}} \Big]+ \\
&+(1\vee \beta)(Q[\mbox{\boldmath$\xi$}]+\|\mbox{\boldmath$\xi$}\|_{\mathfrak{H}}^2+
\|\dot{\mbox{\boldmath$\xi$}}\|_{\mathfrak{H}}^2) = \\
&\geq -\mathfrak{Re}[(\dot{\bm{\xi}}|\bm{\xi})_{\mathfrak{H}}]-\beta \|\dot{\bm{\xi}}\|_{\mathfrak{H}}^2
+(1\vee \beta)\|\bm{\xi}\|_{\mathfrak{H}}^2+
(1\vee \beta)\|\bm{\xi}\|_{\mathfrak{H}}^2 \\
&\geq \|\dot{\bm{\xi}}\|_{\mathfrak{H}}^2
-\mathfrak{Re}[(\dot{\bm{\xi}}|\bm{\xi})_{\mathfrak{H}}]+
\|\bm{\xi}\|_{\mathfrak{H}}^2 \\
&\geq 0,
\end{align*}
since $(\bm{L}\bm{\xi}|\dot{\bm{\xi}})=Q(\bm{\xi}, \dot{\bm{\xi}})=
Q(\dot{\bm{\xi}}, \bm{\xi})^*$ and 
$|\mathfrak{Re}[(\bm{B}\dot{\bm{\xi}}|\dot{\bm{\xi}})_{\mathfrak{H}}| \leq \beta 
\|\dot{\bm{\xi}}\|_{\mathfrak{H}}^2$.

If $\Lambda > \beta$, then the operator $\mathbf{A}+\Lambda$ has the bounded inverse defined on $\mathfrak{E}$. Actually the equation
$$\mathbf{A}U+\Lambda U=F=
\begin{bmatrix}
\mbox{\boldmath$f$} \\
\\
\mbox{\boldmath$g$}
\end{bmatrix}
\in \mathfrak{E}
$$
means
$$
\begin{cases}
-\dot{\mbox{\boldmath$\xi$}}+\Lambda\mbox{\boldmath$\xi$}=\mbox{\boldmath$f$} \in \mathfrak{H}_0^{\mathrm{div}} \\
\\
\bm{B}
\dot{\mbox{\boldmath$\xi$}}+\mbox{\boldmath$L$}\mbox{\boldmath$\xi$}+\Lambda\dot{\mbox{\boldmath$\xi$}}=\mbox{\boldmath$g$} \in \mathfrak{H},
\end{cases}
$$
which can be solved as
$$
\begin{cases}
\mbox{\boldmath$\xi$}=
(\mbox{\boldmath$L$}+\Lambda
\bm{B}
+\Lambda^2)^{-1}(\bm{B}\mbox{\boldmath$f$}+\Lambda\mbox{\boldmath$f$}+\mbox{\boldmath$g$})
\in \mathsf{D}(\mbox{\boldmath$L$}),\\
\\
\dot{\mbox{\boldmath$\xi$}}=
(\mbox{\boldmath$L$}+\Lambda
\bm{B}
+\Lambda^2)^{-1}(\bm{B}\mbox{\boldmath$f$}+\Lambda\mbox{\boldmath$f$}+\mbox{\boldmath$g$})
-\mbox{\boldmath$f$} \in \mathfrak{H}_0^{\mathrm{div}},
\end{cases}
$$
thanks to Proposition \ref{Resolv},
since $\Lambda^2 >|\Lambda\beta|$ for
$\Lambda > \beta$.  $\square$\\

Therefore,
considering the problem
\begin{align}
&\frac{\partial^2\mbox{\boldmath$\xi$}}{\partial t^2}+\bm{B}\frac{\partial\mbox{\boldmath$\xi$}}{\partial t}+
\mbox{\boldmath$L$}\mbox{\boldmath$\xi$}=\bf{0},
\nonumber \\
&\mbox{\boldmath$\xi$}(t,\cdot)\in \mathsf{D}(\mbox{\boldmath$L$})\quad\mbox{for}\quad \forall t \geq 0, \nonumber \\
&\mbox{\boldmath$\xi$}=\overset{\circ}{\bm{\xi}},\quad \frac{\partial\mbox{\boldmath$\xi$}}{\partial t}=\overset{\circ}{\mbox{\boldmath$v$}}
\quad \mbox{at}\quad t=0. \label{HomIBV},
\end{align}
 we can claim

\begin{Theorem}
Suppose $\overset{\circ}{\bm{\xi}} \in \mathsf{D}(\bm{L})$ and 
$\overset{\circ}{\mbox{\boldmath$v$}}\in \mathfrak{H}_0^{\mathrm{div}}$.
 Then the initial-boundary value problem \eqref{HomIBV} admits a unique solution
$$\mbox{\boldmath$\xi$}\in C^2([0,+\infty[, \mathfrak{H})\cap
C^1([0,+\infty[,\mathfrak{H}_0^{\mathrm{div}})
\cap C([0,+\infty[, \mathsf{D}(\mbox{\boldmath$L$}))$$
and the energy
\begin{align}
E(t)&=\|\mbox{\boldmath$\xi$}\|_{\mathfrak{H}^{\mathrm{div}}}^2+\|\dot{\mbox{\boldmath$\xi$}}\|_{\mathfrak{H}}^2 \nonumber \\
&=\|\mbox{\boldmath$\xi$}||_{\mathfrak{H}}^2+Q[\mbox{\boldmath$\xi$}]+\Big\|\frac{\partial\mbox{\boldmath$\xi$}}{\partial t}\Big\|_{\mathfrak{H}}^2
\end{align}
enjoys the estimate
\begin{equation}
\sqrt{E(t)}\leq e^{\Lambda t}\cdot\sqrt{E(0)},
\end{equation}
where $\Lambda=1\vee \beta$.
\end{Theorem}

Here $\mathsf{D}(\mbox{\boldmath$L$})$ is equipped with the norm
$(\|\mbox{\boldmath$\xi$}\|_{\mathfrak{H}^{\mathrm{div}}}^2+\|\mbox{\boldmath$L$}\mbox{\boldmath$\xi$}\|_{\mathfrak{H}}^2)^{1/2}$.\\

Correspondingly we may consider the inhomogeneous initial-boundary value problem
\begin{equation}
\frac{dU}{dt} +\mathbf{A}U=F(t),\quad U|_{t=0}=U_0.  \label{0336}
\end{equation} We can claim

\begin{Proposition}
If $U_0 \in \mathsf{D}(\mathbf{A})$ and $F \in C([0,+\infty[; \mathfrak{E})$, then there exists a unique solution 
$$U\in C^1([0,+\infty[; \mathfrak{E}) \cap C([0,+\infty[; \mathsf{D}(\mathbf{A}) )$$
to the problem \eqref{0336}, and it enjoys the estimate
\begin{equation}
\|U(t)\|_{\mathfrak{E}}\leq 
e^{\Lambda t}\Big(
\|U_0\|_{\mathfrak{E}}+\int_0^t
e^{-\Lambda s}\|F(s)\|_{\mathfrak{E}}ds\Big),
\end{equation}
where $\Lambda=1\vee 2|\Omega|$. 
\end{Proposition}

Therefore,
considering the problem
\begin{align}
&\frac{\partial^2\mbox{\boldmath$\xi$}}{\partial t^2}+\bm{B}\frac{\partial\mbox{\boldmath$\xi$}}{\partial t}+
\mbox{\boldmath$L$}\mbox{\boldmath$\xi$}=\mbox{\boldmath$f$}(t, \mbox{\boldmath$x$}), \nonumber \\
&\mbox{\boldmath$\xi$}=\overset{\circ}{\bm{\xi}},\quad \frac{\partial\mbox{\boldmath$\xi$}}{\partial t}=\overset{\circ}{\mbox{\boldmath$v$}}
\quad \mbox{at}\quad t=0, \nonumber \\
&\mbox{\boldmath$\xi$}(t,\cdot)\in \mathsf{D}(\mbox{\boldmath$L$})\quad\mbox{for}\quad \forall t \geq 0. \label{nonHIBV},
\end{align}
 we can claim

\begin{Theorem}
Suppose $\overset{\circ}{\bm{\xi}} \in \mathsf{D}(\bm{L})$, 
$\overset{\circ}{\mbox{\boldmath$v$}}\in \mathfrak{H}_0^{\mathrm{div}}$
and $ \mbox{\boldmath$f$} \in C([0,+\infty[; \mathfrak{H})$. Then the initial-boundary value problem \eqref{nonHIBV} admits a unique solution
$$\mbox{\boldmath$\xi$}\in C^2([0,+\infty[, \mathfrak{H})\cap
C^1([0,+\infty[,\mathfrak{H}_0^{\mathrm{div}})
\cap C([0,+\infty[, \mathsf{D}(\mbox{\boldmath$L$}))$$
and the energy
\begin{align}
E(t)=E(t, \mbox{\boldmath$\xi$})&:=\|\mbox{\boldmath$\xi$}\|_{\mathfrak{H}^{\mathrm{div}}}^2+\|\dot{\mbox{\boldmath$\xi$}}\|_{\mathfrak{H}}^2 \nonumber \\
&=\|\mbox{\boldmath$\xi$}||_{\mathfrak{H}}^2+Q[\mbox{\boldmath$\xi$}]+\Big\|\frac{\partial\mbox{\boldmath$\xi$}}{\partial t}\Big\|_{\mathfrak{H}}^2
\end{align}
enjoys the estimate
\begin{equation}
\sqrt{E(t)}\leq e^{\Lambda t}\Big(\sqrt{E(0)}+\int_0^te^{-\Lambda s}
\|\mbox{\boldmath$f$}(s)\|_{\mathfrak{H}}ds \Big)
\end{equation}
for $\Lambda=1\vee \beta$.
\end{Theorem}



\section{Eigenfrequency, eigenvector, the variational principle}

Astrophysicists used to discuss on the so called `variational principle'. See \cite{Chandra64}, \cite{Clement}, \cite{LyndenBO}, and so on. Although they discuss
about self-gravitating gaseous masses, we would like to follow their discussions by applying them to the present case of the model of rotating atmosphere on the Earth, namely, we consider the linearized wave equation \eqref{LEAMa}.\\

Since the operator 
$\displaystyle \bm{J}: \bm{v} \mapsto 
\begin{bmatrix}
0 & -1 & 0 \\
1 & 0 & 0\\
0 & 0 & 0
\end{bmatrix}\bm{v}$ is skew symmetric, we introduce the operator $\mbox{\boldmath$J$}_{\star}$ defined by
\begin{equation}
\mbox{\boldmath$J$}_{\star}=\mathrm{i}\mbox{\boldmath$J$}.
\end{equation}
Then the operator $\mbox{\boldmath$J$}_{\star}$ is a bounded self-adjoint operator on $\mathfrak{H}$. 
Recall 
\begin{equation}
(\mbox{\boldmath$J$}_{\star}\mbox{\boldmath$\xi$}|\mbox{\boldmath$\xi$})_{\mathfrak{H}}=
-2\int_{\mathfrak{R}}\mathfrak{Im}[\xi^1(\xi^2)^*]\bar{\rho}(\mbox{\boldmath$x$})d\mbox{\boldmath$x$}. \label{V.2}
\end{equation}
The equation
\eqref{LEAMa}
reads
\begin{equation}
\frac{\partial^2\mbox{\boldmath$\xi$}}{\partial t^2}-2\Omega\mathrm{i}\mbox{\boldmath$J$}_{\star}\frac{\partial\mbox{\boldmath$\xi$}}{\partial t}
+\mbox{\boldmath$L$}\mbox{\boldmath$\xi$}=0 \label{V.3}
\end{equation}\\

Let us suppose that there exists a solution ${\mbox{\boldmath$\xi$}}$ to \eqref{V.3} of the form
\begin{equation}
{\mbox{\boldmath$\xi$}}(t,\mbox{\boldmath$x$})=e^{\mathrm{i}\sigma t}\mbox{\boldmath$\Xi$}(\mbox{\boldmath$x$}) \label{V.4}
\end{equation}
where $\sigma \in \mathbb{C}$ and $\mbox{\boldmath$\Xi$}\in \mathsf{D}({\mbox{\boldmath$L$}})$. Then the equation \eqref{V.3} reduces to
\begin{equation}
-\sigma^2\mbox{\boldmath$\Xi$}+\sigma 2\Omega \mbox{\boldmath$J$}_{\star}\mbox{\boldmath$\Xi$}+\mbox{\boldmath$L$}\mbox{\boldmath$\Xi$}=\mathbf{0}, \label{EigEq}
\end{equation}
or the equation \eqref{LEAMa}reduces to
\begin{equation}
-\sigma^2\mbox{\boldmath$\Xi$}+\mathrm{i}\sigma 2\Omega \mbox{\boldmath$J$}\mbox{\boldmath$\Xi$}+\mbox{\boldmath$L$}\mbox{\boldmath$\Xi$}=\mathbf{0}, \label{non*EigEq}.
\end{equation}

So, we use 

\begin{Definition}
When the equation \eqref{EigEq} [\!( \eqref{non*EigEq} )\!] is satisfied for $\sigma \in \mathbb{C} $ and $\bm{\Xi} \in \mathsf{D}(\bm{L}), \not=\mathbf{0}$, then $\sigma$ is called an eigenfrequency of the wave equation 
\eqref{V.3} [\!( \eqref{LEAMa})\!] and $\bm{\Xi}$ is called an eigenvector associated with the eigenfrequency $\sigma$.
\end{Definition}

Note that $0$ is an eigenfrequency. In fact the vector
$$\bm{\Xi}=\frac{1}{\bar{\rho}}\nabla \times \bm{a} $$
belongs to $\mathrm{Ker}\bm{L}$ for any $\bm{a} \in C_0^{\infty}(\mathfrak{R}; \mathbb{C}^3) $, and turns out to be an eigenvector associated with the eingenfrequency $0$ provided that $\bm{\Xi}\not=\mathbf{0}$.\\

Here let us recall the operator $\mathbf{A}$ defined as
$$
\mathbf{A}U=
\begin{bmatrix}
O & -I \\
\\
\bm{L} & 2\Omega\bm{J}
\end{bmatrix}
U=
\begin{bmatrix}
-\dot{\bm{\Xi}} \\
\\
2\Omega\bm{J}\dot{\bm{\Xi}}+\bm{L}\bm{\Xi}
\end{bmatrix}
$$
for
$$
U=\begin{bmatrix}
\bm{\Xi} \\
\\
\dot{\bm{\Xi}}
\end{bmatrix}
\in \mathsf{D}(\mathbf{A})=\mathsf{D}(\bm{L})\times\frak{H}_0^{\mathrm{div}},
$$
which was introduced in order to alpply the Hille-Yosida theory to the linear evolution equation in the preceeding section. Obviously we can claim

\begin{Proposition}
If and only if $\sigma \in \mathbb{C}$ is an eigenfrequency of the equation \eqref{V.3} [\!( \eqref{LEAMa})\!] ,
$\Lambda=\mathrm{i}\sigma $ is an eigenvalue of the operator $-\mathbf{A}$.
\end{Proposition}

Let us introduce the following

\begin{Definition}
We denote by $\mathfrak{L}$ the one parameter family of operators
$(-\sigma^2+\sigma 2\Omega\bm{J}_{\star}+\bm{L})_{\sigma \in \mathbb{C}}$, which is called 'quadratic pencil'. We denote 
\begin{equation}
\mathfrak{L}(\sigma):=
-\sigma^2+\sigma 2\Omega\bm{J}_{\star}+\bm{L}.
\end{equation}

If the operator $\mathfrak{L}(\sigma)=-\sigma^2+\sigma 2\Omega \bm{J}_{\star}+\bm{L}$ admits the bounded inverse defined on $\mathfrak{H}$, $\sigma$ is said to belong to the resolvent set $\varrho(\mathfrak{L})$. We denote
$\sigma(\mathfrak{L})=\mathbb{C} \setminus \varrho(\mathfrak{L})$, and call it the spectrum of the quadratic pencil $\mathfrak{L}$. 
\end{Definition}

If $\sigma$ is an eigenfrequency, then it belongs to the spectrum $\sigma(\mathfrak{L})$, so, $ 0 \in \sigma(\mathfrak{L})$, but belonging to
$\sigma(\mathfrak{L})$ does not mean being an eigenfrequency a priori, of course. 

\begin{Proposition}
$\varrho(\mathfrak{L})$ is an open subset of $\mathbb{C}$, and $\sigma(\mathfrak{L})$ is closed.
\end{Proposition}

Proof. Let us consider $\sigma \in \varrho(\mathfrak{L})$. Then 
$$\mathfrak{L}(\sigma+\Delta\sigma)=
\mathfrak{L}(\sigma)\Big[
I+ \mathfrak{L}(\sigma)^{-1}(\Delta\sigma(-2\sigma+2\Omega\bm{J}_{\star}-\Delta\sigma)\Big]
$$
admits the bounded inverse and $\sigma+\Delta\sigma \in \varrho(\mathfrak{L})$, if
$$|\|\mathfrak{L}(\sigma)^{-1}\Delta\sigma(-2\sigma+2\Omega\bm{J}_{\star}-\Delta\sigma)\||_{\mathcal{B}(\mathfrak{H})} <1.$$
For this inequality, it is sufficient that
$$|\|\mathfrak{L}(\sigma)^{-1}\||_{\mathcal{B}(\mathfrak{H})}\cdot|\Delta\sigma|\cdot|(2(|\sigma|
+|\Omega|)+|\Delta\sigma| ) < 1,$$
or
$$|\Delta\sigma|<-(|\sigma|+|\Omega|)+\sqrt{\frac{1}{|\|\mathfrak{L}(\sigma)^{-1}\||_{\mathcal{B}(\mathfrak{H})}}+
(|\sigma|+|\Omega|)^2}.$$
This means $\varrho(\mathfrak{L})$ is open.
$\square$\\

We claim

\begin{Proposition}
It holds that
\begin{equation}
\mathrm{i}\varrho(\mathfrak{L})=\varrho(-\mathbf{A}),
\quad
\mathrm{i}\sigma(\mathfrak{L})=\sigma(-\mathbf{A}).
\end{equation}
Here $\varrho(-\mathbf{A}), \sigma(-\mathbf{A})$ stand for the usual resolvent set, the spectrum of the operator $-\mathbf{A}$ in the Hilbert space $\mathfrak{E}=\mathfrak{H}_0^{\mathrm{div}}\times \mathfrak{H}$.
\end{Proposition}

Proof. Let $\sigma \in \varrho(\mathfrak{L})$ and $\Lambda=\mathrm{i}\sigma$. Consider the equation
$$\mathbf{A}U+\Lambda U=F=
\begin{bmatrix}
\bm{f} \\
\\
\bm{g}
\end{bmatrix}
\in \mathfrak{E}, $$
or,
$$
\begin{cases}
&-\dot{\bm{\Xi}}+\mathrm{i}\sigma \bm{\Xi}=\bm{f} \in \mathfrak{H}_0^{\mathrm{div}}, \\
& \\
&2\Omega J\dot{\bm{\Xi}}+\bm{L}\bm{\Xi}+
\mathrm{i}\sigma\dot{\bm{\Xi}}=\bm{g} \in \mathfrak{H}.
\end{cases}
$$
This system of equations can be solved
as
$$
\begin{cases}
\bm{\Xi}&=
(-\sigma^2-\mathrm{i}\sigma 2\Omega\bm{J}+\bm{L})^{-1}
(2\Omega J \bm{f}-\mathrm{i}\sigma\bm{f}+\bm{g})
\in \mathsf{D}(\bm{L}), \\
& \\
\dot{\bm{\Xi}}&=
(-\sigma^2-\mathrm{i}\sigma 2\Omega\bm{J}+\bm{L})^{-1}
(2\Omega J \bm{f}-\mathrm{i}\sigma\bm{f}+\bm{g})
-\bm{f} \in \mathfrak{H}_0^{\mathrm{div}},
\end{cases}
$$
since $\sigma \in \varrho(\mathfrak{L})$, while
$F \mapsto U$ is continuous. Therefore $\Lambda=\mathrm{i}\sigma \in 
\varrho(-\mathbf{A})$, or, $\mathrm{i}\varrho(\mathfrak{L}) \subset
\varrho(-\mathbf{A})$.

Inversely let $\Lambda \in \varrho(-\mathbf{A})$ and $ \sigma=-\mathrm{i}\Lambda$. Consider the equation
$$(-\sigma^2+\mathrm{i}\sigma 2\Omega\bm{J}+\bm{L})\bm{\Xi}=\bm{f}
\in \mathfrak{H},$$
or
$$(\Lambda^2+\Lambda 2\Omega\bm{J}+\bm{L})\bm{\Xi}=\bm{f},
$$
which is equivalent to the system of equations
$$
\begin{cases}
& \Lambda\dot{\bm{\Xi}}+2\Omega\bm{J}\dot{\bm{\Xi}}
+\bm{L}\bm{\Xi}=\bm{f}, \\
& \\
&\dot{\bm{\Xi}}=\Lambda\bm{\Xi}.
\end{cases}
$$
But this is nothing but
$$\mathbf{A}U+\Lambda U=
\begin{bmatrix}
\mathbf{0} \\
\\
\bm{f}
\end{bmatrix}.
$$
Since $\Lambda \in \varrho(-\mathbf{A})$ is supposed, this admits the solution
$$U=
\begin{bmatrix}
\bm{\Xi} \\
\\
\dot{\bm{\Xi}}
\end{bmatrix}=
(\mathbf{A}+\Lambda)^{-1}
\begin{bmatrix}
\mathbf{0} \\
\\
\bm{f}
\end{bmatrix},
$$
and $\bm{f} \mapsto \bm{\Xi}$ is continuous, that is,
$\sigma =-\mathrm{i}\Lambda \in \varrho(\mathfrak{L})$. $\square$\\

When $\Omega=0$, we have 
$$\sigma(\mathfrak{L})=
\{ \sigma \in \mathbb{C} \   |\   \lambda=\sigma^2 \in \sigma(\bm{L}) \},
$$
where $\sigma(\bm{L})$ is the spectrum of the self-adjoint operator $\bm{L}$
in the Hilbert space $\mathfrak{H}$. Since $\bf{L}$ is self-adjoint and $\bm{L}\geq 0$, we have $\sigma(\bm{L}) \subset \{ \lambda \in \mathbb{R} \  |\  \lambda \geq 0\}$. Thus it holds that
\begin{equation}
\sigma(\mathfrak{L}) \subset \mathbb{R} \quad\mbox{and}\quad
\{\Lambda\in\mathbb{C}\  |\  \mathfrak{Re}[\Lambda]\not=0\}\subset\varrho(-\mathbf{A}),\label{0508}
\end{equation}
when $\Omega=0$.
But, when $\Omega\not=0$, the situation is not so evident. At least, 
we can claim
$$\sigma(\mathfrak{L}) \subset
\mathbb{C}\setminus
]-\infty, -2|\Omega|[\mathrm{i}, $$
since
$$\{ \Lambda \in \mathbb{R}\  |\  \Lambda >2|\Omega| \}
\subset \varrho(-\mathbf{A}), $$
as shown in the proof of Proposition \ref{Prop.2}. 
Moreover, it is clear that  $\mathfrak{L}(\sigma^*)=(\mathfrak{L}(\sigma))^*$, since $2\Omega\bm{J}_{\star}$ and $\bm{L}$ are self-adjoint. Therefore we have
\begin{equation}
\sigma \in \varrho(\mathfrak{L})\quad\Leftrightarrow\quad
\sigma^*\in\varrho(\mathfrak{L}).
\end{equation}
This means that $\varrho(\mathfrak{L})$ and $\sigma(\mathfrak{L})$ are symmetric re the real axis in the complex number plane. Correspondingly, $\sigma(-\mathbf{A})$ and $\varrho(-\mathbf{A})$ are symmetric re the imaginary axis. Thus we can claim
\begin{align}
&\sigma(\mathfrak{L}) \subset
\mathbb{C}\setminus
(]-\infty, -2|\Omega|[
\cup ]2|\Omega|, +\infty[)\mathrm{i} \quad\mbox{and} \nonumber \\
&\{ \Lambda \in \mathbb{R}\  |\  |\Lambda| >2|\Omega| \} \subset
\varrho(-\mathbf{A}). \label{0510}
\end{align}
However 
 the gap between the information \eqref{0508} for $\Omega=0$ and that \eqref{0510} for $\Omega\not=0$ is too much. So, we are want to strengthen \eqref{0510}  when $\Omega\not=0$. In order to do it, we use the following
 
 \begin{Proposition}\label{Rbelow}
 If $\sigma \in \varrho(\mathfrak{L})$, then it holds that
 \begin{equation}
 |\| \mathfrak{L}(\sigma)^{-1}\||_{\mathcal{B}(\mathfrak{H})} \geq \frac{1}{d(2(|\sigma|+|\Omega|)+d)},
 \end{equation}
 where $d:=\mathrm{dist}(\sigma, \sigma(\mathfrak{L}))$.
 \end{Proposition}

Proof. Let $\sigma \in \varrho(\mathfrak{L})$. Then, for $\Delta\sigma \in \mathbb{C}$, the operator
$$\mathfrak{L}(\sigma+\Delta\sigma)=
\mathfrak{L}(\sigma)\Big[I+
\mathfrak{L}(\sigma)^{-1}\Delta\sigma(-2\sigma+2\Omega\bm{J}_{\star}-\Delta\sigma)\Big]$$
admits the bounded inverse in $\mathcal{B}(\mathfrak{H})$ and $\sigma+\Delta\sigma \in
\varrho(\mathfrak{L})$, if
$$|\|\mathfrak{L}(\sigma)^{-1}\Delta\sigma(-2\sigma+2\Omega\bm{J}_{\star}-\Delta\sigma)\||_{\mathcal{B}(\mathfrak{H})} <1.$$
For this inequality, it is sufficient that
$$|\|\mathfrak{L}(\sigma)^{-1}\||_{\mathcal{B}(\mathfrak{H})}\cdot|\Delta\sigma|\cdot|(-2\sigma
+2\Omega\bm{J}_{\star}-\Delta\sigma| <1.$$
In other words, if  $\sigma+\Delta\sigma \in \sigma(\mathfrak{L})$, then it should hold
$$|\|\mathfrak{L}(\sigma)^{-1}\||_{\mathcal{B}(\mathfrak{H})}\cdot|\Delta\sigma|\cdot|(-2\sigma
+2\Omega\bm{J}_{\star}-\Delta\sigma| \geq 1,
$$
and necessarily
$$|\|\mathfrak{L}(\sigma)^{-1}\||_{\mathcal{B}(\mathfrak{H})}\cdot|\Delta\sigma|\cdot|(2(|\sigma|
+|\Omega|)+|\Delta\sigma| )\geq 1.$$
If $d <+\infty$, then there is a sequence $\sigma+(\Delta\sigma)_n
\in \sigma(\mathfrak{L})$ such that $|(\Delta\sigma)_n| \rightarrow d$, and the assertion follows. $\square$\\

Let us fix $\sigma_{\infty} \in \partial\sigma(\mathfrak{L})$. We are going to show $\sigma_{\infty} \in \mathbb{R}$.

Let us consider a sequance $(\sigma_n)_n$ such that $\sigma_n \in \varrho(\mathfrak{L})$ and
$\sigma_n \rightarrow \sigma_{\infty}$ as $ n \rightarrow \infty$. By Proposition \ref{Rbelow} 
we have $|\|(\mathfrak{L}(\sigma_n)^{-1}\||_{\mathcal{B}(\mathfrak{H})} \rightarrow +\infty$, therefore
there are $\bm{f}_n \in \mathfrak{H}$ such that $\|\bm{f}_n\|_{\mathfrak{H}}=1$ and $\|\mathfrak{L}(\sigma_n))^{-1}\bm{f}_n\|_{\mathfrak{H}}
 \rightarrow +\infty$ as $n \rightarrow \infty$. Put $\bm{\xi}_n=\mathfrak{L}(\sigma_n)^{-1}\bm{f}_n (\in \mathsf{D}(\bm{L}))$ and $\bm{\eta}_n=\bm{\xi}_n/\|\bm{\xi}_n\|_{\mathfrak{H}}$. Then $\|\bm{\eta}_n\|_{\mathfrak{H}}=1$ and $$
 \Big|(\mathfrak{L}(\sigma_n)\bm{\eta}_n|\bm{\eta}_n)_{\mathfrak{H}}\Big|
 =\Big|\frac{1}{\|\bm{\xi}_n\|_{\mathfrak{H}}^2}(\bm{f}_n|\bm{\xi}_n)_{\mathfrak{H}}\Big|
 \leq \frac{1}{\|\bm{\xi}_n\|_{\mathfrak{H}}} \rightarrow 0.$$
 But we see
 $$
 (\mathfrak{L}(\sigma_n)\bm{\eta}_n|\bm{\eta}_n)_{\mathfrak{H}}=-(\sigma_n)^2+\sigma_nb_n+c_n,$$
 where
 $$b_n:=2\Omega(\bm{J}_{\star}\bm{\eta}_n|\bm{\eta}_n)_{\mathfrak{H}},
 \quad c_n=(\bm{L}\bm{\eta}_n|\bm{\eta}_n)_{\mathfrak{H}}=Q[\bm{\eta}_n].$$
 Therefore $b_n, c_n $ are real and 
 $$|b_n|\leq 2|\Omega|,\qquad c_n \geq 0.$$
 Hence, by taking a subsequence if necessary, we can suppose that $b_n$ tends to a limit $b_{\infty}$ such that $|b_{\infty}|\leq2|\Omega|$. Put
 $c_{\infty}:=(\sigma_{\infty})^2-\sigma_{\infty}b_{\infty}$. Then we see 
 $c_n \rightarrow c_{\infty}$. Hence $c_{\infty}$ is real and $ \geq 0$,
 and $\sigma_{\infty}$ turns out to enjoy the quadratic equation
 $$-(\sigma_{\infty})^2+b_{\infty}\sigma_{\infty}+c_{\infty}=0.$$
 Consequently, 
 $$\sigma_{\infty}=\frac{b_{\infty}}{2}+\sqrt{\frac{b_{\infty}^2}{4}+c_\infty}
 \quad\mbox{or}\quad
 \sigma_{\infty}=\frac{b_{\infty}}{2}-\sqrt{\frac{b_{\infty}^2}{4}+c_\infty},$$
 so, anyway, $\sigma_{\infty} \in \mathbb{R}$. This was to be prooved. \\
 
 Summing up, we can claim
 \begin{Proposition}\label{BoundarySpectrum}
 It holds that
 \begin{equation}
\partial  \sigma(\mathfrak{L}) \subset \mathbb{R}. \label{BdSpec}
 \end{equation}
 \end{Proposition}
 
 This conclusion owes to \cite[Theorem 1]{DysonS}. But their original proof 
 is little bit logically weak, and we have needed to edit it as above. \\
 
 Let us consider $\sigma_0=\alpha_0+\beta_0\mathrm{i} \in \sigma(\mathfrak{L})$, where
 $\alpha_0, \beta_0 \in \mathbb{R}$. We are going to show that $|\beta_0|\not=0$ implies a contradiction. By the symmetricity of $\sigma(\mathfrak{L})$ we can suppose $\beta_0 >0$ without loss of generality. Choosing $K >2|\Omega| \vee \beta_0$, we consider the segment
 \begin{align*}
 I&=[\sigma, K\mathrm{i}] \\
 &=\{ \sigma(t)=(1-t)\alpha_0+
 (\beta_0+(K-\beta_0)t)\mathrm{i} \  |\  0\leq t \leq 1 \}.
 \end{align*}
 Note that $\sigma(0)=\sigma_0 \in \sigma(\mathfrak{L})$ and
 $\sigma(1)=K\mathrm{i} \in \varrho(\mathfrak{L})$ by \eqref{0510}, since 
 $K >2|\Omega|$. Put
 $$\bar{t}:=\sup\{ t\in [0,1]\  |\  \sigma(t) \in \sigma(\mathfrak{L}) \}.$$
 Then $0\leq \bar{t} <1, \sigma(\bar{t})\in \sigma(\mathfrak{L})$ and
 $\sigma(t) \in \varrho(\mathfrak{L})$ for $t >\bar{t}$. Hence $\sigma(\bar{t})\in \partial\sigma(\mathfrak{L})$. But $\mathfrak{Im}[\sigma(\bar{t})] \geq \beta_0>0$, a contradiction to $\partial\sigma(\mathfrak{L}) \in \mathbb{R}$. $\square$\\
 
 Therefore we can claim

\begin{Theorem} \label{Th.realspector}
It holds that
\begin{equation}
\sigma(\mathfrak{L}) \subset \mathbb{R}
\end{equation}
even when $\Omega\not=0$.
\end{Theorem}


Now we note that it is known that there is a sequence of eigenfrequencies $\sigma_{n, \pm}=
\pm \sqrt{\lambda_n^{\mathcal{N}}}, n \in \mathbb{N}$, when $\Omega=0$. ( See the discussion given later.) However, up to now, we have no knowledge on the existence of eigenfrequencies when $\Omega\not=0$.
When $\Omega=0$, then \eqref{non*EigEq} reads
$$ -\sigma^2\bm{\Xi}+\bm{L}\bm{\Xi}=\mathbf{0}, $$
so that the eigenvctor $\bm{\Xi}$ associated with the eigenfrequency $\sigma\not=0$ can be supposed to be real, since $\sigma \in \mathbb{R}$. However, when $\Omega\not=0$, the situation is different. Namely, we claim

\begin{Proposition}
Suppose $\Omega\not=0$. Let $\sigma\not=0$ be an eigenfrequency of
the equation \eqref{V.3} [\!( \eqref{LEAMa})\!] and $\bm{\Xi}$ be an associated eigenvector. Then the eigenvector  $\bm{\Xi}$ is impossible to be real, that is, $\mathfrak{Im}[\bm{\Xi}(\bm{x})]$ cannot vanish everywhere. 
\end{Proposition}

Proof. By Theorem \ref{Th.realspector} we see $\sigma \in \mathbb{R}$. Let us look at \eqref{non*EigEq}:
\begin{equation}
-\sigma^2\bm{\Xi}+\mathrm{i}\sigma 2\bm{\Omega}\times\bm{\Xi}+\bm{L}\bm{\Xi}=\mathbf{0}. \label{R1}.
\end{equation}
Let us denote $\bm{X}(\bm{x})=\mathfrak{Re}[\bm{\Xi}(\bm{x})],
\bm{Y}(\bm{x})=\mathfrak{Im}[\bm{\Xi}(\bm{x})]$ so that
$\bm{X}(\bm{x}), \bm{Y}(\bm{x}) \in \mathbb{R}, \bm{\Xi}=\bm{X}+\mathrm{i}\bm{Y}$, 
Suppose that $\bm{Y}=0$ and deduce a contradiction. Now \eqref{R1} means
\begin{subequations}
\begin{align}
&-\sigma^2\bm{X}-\sigma 2\bm{\Omega}\times \bm{Y}+\bm{L}\bm{X}=\mathbf{0}, \label{R2a} \\
&-\sigma\bm{Y}+\sigma 2\bm{\Omega}\times \bm{X}+\bm{L}\bm{Y}=\mathbf{0}. \label{R2b}
\end{align}
\end{subequations}
Since $\bm{Y}=\mathbf{0}$ is supposed, this reads
\begin{subequations}
\begin{align}
&-\sigma^2\bm{X}+\bm{L}\bm{X}=\mathbf{0}, \label{R3a} \\
&\sigma 2\bm{\Omega}\times \bm{X}=\mathbf{0} \label{R3b}
\end{align}
\end{subequations}
Since $\sigma\not=0, \Omega\not=0$, \eqref{R3b} implies
\begin{equation}
X^1=X^2=0, \label{R4}
\end{equation}
where $\bm{X}=(X^1, X^2, X^3)^{\top}$. Then \eqref{R3a} reads
\begin{subequations}
\begin{align}
\frac{\partial}{\partial x^1}\Big(-\overline{\frac{d\varUpsilon}{d\rho}}
\frac{\partial}{\partial x^3}(\bar{\rho}X^3)\Big)&=0, \label{R5a} \\
 \frac{\partial}{\partial x^2}\Big(-\overline{\frac{d\varUpsilon}{d\rho}}
\frac{\partial}{\partial x^3}(\bar{\rho}X^3)\Big)&=0, \label{R5b} \\
-\sigma^2 X^3 +
\frac{\partial}{\partial x^3}\Big(-\overline{\frac{d\varUpsilon}{d\rho}}
\frac{\partial}{\partial x^3}(\bar{\rho}X^3)\Big)&=0. \label{R5c} 
\end{align}
\end{subequations}
Consequently, \eqref{R5a} and \eqref{R5b} imply that 
$\displaystyle -\overline{\frac{d\varUpsilon}{d\rho}}\frac{\partial}{\partial x^3}(\bar{\rho}X^3)$ is a function of $x^3$ independent of
$x^1, x^2$, and,
since $\sigma\not=0$, \eqref{R5c} implies that 
$X^3$ is so, too.
 However $\bm{\Xi}=\bm{X} \in \mathsf{D}(\bm{L})$ suppose the boundary condition
$$\Big(\bm{\Xi}\Big|\frac{\partial}{\partial r}\Big)=0\quad \mbox{on}\quad r=R_0.$$
Namely, 
$$X^3(x^3)x^3=0\quad\mbox{on}\quad \Sigma_0=\{\   r= R_0\  \},
$$
therefore $X^3=0$, and $\bm{\Xi}=(X^1,X^2.X^3)^{\top}=\mathbf{0}$, a contradiction. $\square$\\

Let $\sigma$ be an eigenfrequency of the equation \eqref{LEAMa}and $\bm{\Xi}$ be an associated eigenvector. 

Then $\bm{\xi}(t,\bm{x})=e^{\mathrm{i}\sigma t}\bm{\Xi}(\bm{x})$ is a solution of the equation \eqref{3.12}. Since the coefficients of the equation \eqref{LEAMa}are real and the equation is linear, we can claim that $\bm{\xi}(t,\bm{x})^*, \mathfrak{Re}[\bm{\xi}(t,\bm{x})],
\mathfrak{Im}[\bm{\xi}(t,\bm{x})]$ are solutions of \eqref{3.12}, too. But, since $\sigma$ is real by Theorem \ref{Th.realspector}, we see
\begin{align*}
\bm{\xi}_{\mathfrak{R}}(t,\bm{x}):=\mathfrak{Re}[\bm{\xi}(t,\bm{x})]&=\mathfrak{Re}[e^{\mathrm{i}\sigma t}\bm{\Xi}(\bm{x})]= \\
&=\cos(\sigma t) \bm{\Xi}_{\mathfrak{R}}(\bm{x})-
\sin(\sigma t)\bm{\Xi}_{\mathfrak{I}}(\bm{x}) = \\
&=
\begin{bmatrix}
|\Xi(\bm{x})^1|\cos(\sigma t+\alpha_1(\bm{x})) \\
\\
|\Xi(\bm{x})^2|\cos(\sigma t+\alpha_2(\bm{x})) \\
\\
|\Xi(\bm{x})^3|\cos(\sigma t+\alpha_3(\bm{x})) 
\end{bmatrix},
\end{align*}
where
\begin{align*}
& \bm{\Xi}_{\mathfrak{R}}(\bm{x}):=\mathfrak{Re}[\bm{\Xi}(\bm{x})],\quad
\bm{\Xi}_{\mathfrak{I}}(\bm{x}):=\mathfrak{Im}[\bm{\Xi}(\bm{x})] \\
& \tan \alpha_j(\bm{x})=\frac{{\Xi}_{\mathfrak{I}}(\bm{x})^j}{{\Xi}_{\mathfrak{R}}(\bm{x})^j},\quad\mbox{so that}\quad 
{\Xi}(\bm{x})^j=|{\Xi}(\bm{x})^j|e^{\mathrm{i}\alpha_j(\bm{x})}.
\end{align*}

Note that the field $\bm{\xi}_{\mathfrak{R}}(t,\bm{x})$ is a real-valued solution of \eqref{LEAMa}
such that $\bm{\xi}_{\mathfrak{R}}(0,\bm{x})=\bm{\Xi}_{\mathfrak{R}}(\bm{x})$.\\

Let $\sigma, \mbox{\boldmath$\Xi$}(\not=\mathbf{0})$ be an eigenfrequency and an associated eigenvector.
Multiplying \eqref{EigEq} by $\mbox{\boldmath$\Xi$}^*\bar{\rho}(\mbox{\boldmath$x$})$ and integrating it, we have
\begin{equation}
-\sigma^2\|\mbox{\boldmath$\Xi$}\|_{\mathfrak{H}}^2+
\sigma 2\Omega(\mbox{\boldmath$J$}_{\star}\mbox{\boldmath$\Xi$}|\mbox{\boldmath$\Xi$})_{\mathfrak{H}}+Q[\mbox{\boldmath$\Xi$}]=0. \label{V.6}
\end{equation} 
Recall the quadratic form  $Q$ is defined by \eqref{3.19}. If we write
\begin{equation}
a=\|\mbox{\boldmath$\Xi$}\|_{\mathfrak{H}}^2,\quad
b=2\Omega (\mbox{\boldmath$J$}_{\star}\mbox{\boldmath$\Xi$}|\mbox{\boldmath$\Xi$})_{\mathfrak{H}},
\quad c=Q[\mbox{\boldmath$\Xi$}], \label{V.7}
\end{equation}
then $a, b, c$ are real numbers, and $\sigma$ satisfies
the quadratic equation
\begin{equation}
-a\sigma^2+b\sigma +c=0, \label{V.8}
\end{equation}
whoe roots are
\begin{equation}
\sigma=\frac{b}{2a}\pm\sqrt{\frac{b^2}{4a^2}+\frac{c}{a}}.\label{V.9}
\end{equation}\\

Here D. Lynden-Bell and J. P. Ostriker, \cite[p.301, line 18]{LyndenBO}, say:
\begin{quote}
Equation (36) [ read \eqref{V.9} ] shows that the system is stable if $c$ is positive for each eigen $\mbox{\boldmath$\xi$}$ [ read $\mbox{\boldmath$\Xi$}$ ] . This assured if $\mathbf{C}$ [ read $\mbox{\boldmath$L$}$] is positive definite. Thus:

A sufficient condition for stability is that $\mathbf{C}$ [read $\mbox{\boldmath$L$}$ ] is positive definite. This is {\it the } condition for secular stability.
\end{quote}

This saying sounds strange. In fact, we may suppose the meaning of the words `stability'
and `secular stability' as C. Hunter \cite{Hunter} defines:

\begin{quote}
A general system is said to be ordinarily or dynamically unstable if the amplitude of some mode grows exponentially in time, but ordinarily stable if every mode is oscillatory in time. An ordinarily stable system can be said to be secularly unstable if small additional dissipative forces can cause some perturbation to grow. Otherwise, the system is sad to be secularly stable.
\end{quote}

As C. Hunter says in the same article, the definition of secular instability does not always confirm to that given above, so, we consider the (ordinary) stability. If there is an eigenfrequancy $\sigma$ which is not real, then the system described by \eqref{V.3} is unstable. It is true. But this means that there is an eigenfrequency $\sigma$ and an associated eigenvector $\mbox{\boldmath$\Xi$}$ such that
$$\frac{b^2}{4a^2}+\frac{c}{a} <0.$$
Normalizing $a=\|\mbox{\boldmath$\Xi$}\|_{\mathfrak{H}}^2=1$, this means 
$$\frac{b^2}{4}+c <0.$$
Therefore we can claim that the system is unstable if there is an eigenvector $\mbox{\boldmath$\Xi$}$ such that
$$\frac{b^2}{4}+c =\Omega^2
((\mbox{\boldmath$J$}_{\star}\mbox{\boldmath$\Xi$}|\mbox{\boldmath$\Xi$})_{\mathfrak{H}}
)^2
+Q[\mbox{\boldmath$\Xi$}] <0.
$$
Of course, in the situation under consideration, this cannot happen, since $Q[\bm{\Xi}]\geq 0$ for $\forall \bm{\Xi} \in \mathsf{D}(\bm{L})$ and $b \in \mathbb{R}$.

However, logically speaking,  the condition $
c=Q[\mbox{\boldmath$\Xi$}]>0 $  for each eigenvector $\mbox{\boldmath$\Xi$}$ is far from the condition of the stability,
contrary to the saying of D. Lynden-Bell and J. P. Ostriker. 

Moreover let us note the following fact:  $Q[\mbox{\boldmath$\Xi$}] \geq 0$ for any
$\mbox{\boldmath$\Xi$} \in \mathfrak{H}^{\mathrm{div}}$, but $Q[\mbox{\boldmath$\Xi$}]=0$ does not imply $\mbox{\boldmath$\Xi$}=\mathbf{0}$; In fact 
$$\mbox{\boldmath$\Xi$}(\mbox{\boldmath$x$})=\frac{1}{\bar{\rho}(\mbox{\boldmath$x$})}\nabla\times \mbox{\boldmath$a$},$$
$\mbox{\boldmath$a$}$ being an arbitrary vector field belonging to
$C_0^{\infty}(\mathfrak{R})$, belongs to the kernel of $\mbox{\boldmath$L$}$, that is, $Q[\mbox{\boldmath$\Xi$}]=0$.\\

Anyway, we are going to describe the `variational principle'.

Let us suppose that there exist $\sigma_0 \in \mathbb{R}$ and $\mbox{\boldmath$\Xi$}_0 \in \mathfrak{H}^{\mathrm{div}}\ (\not=\mathbf{0})$ such that
\begin{equation}
-\sigma_0^2\|\mbox{\boldmath$\Xi$}_0\|_{\mathfrak{H}}^2+
\sigma_0 2\Omega (\mbox{\boldmath$J$}_{\star}\mbox{\boldmath$\Xi$}_0|\mbox{\boldmath$\Xi$}_0)_{\mathfrak{H}}+
Q[\mbox{\boldmath$\Xi$}_0]=0. \label{V.10}
\end{equation}
Of course if $\sigma_0, \mbox{\boldmath$\Xi$}_0$ are real eigenfrequency and an associated eigenvector then \eqref{V.10} is satisfied. Now we assume that
\begin{equation}
a_0=\|\mbox{\boldmath$\Xi$}_0\|_{\mathfrak{H}}^2,\quad
b_0=2\Omega
 (\mbox{\boldmath$J$}_{\star}\mbox{\boldmath$\Xi$}_0|\mbox{\boldmath$\Xi$}_0)_{\mathfrak{H}},\quad
c_0=Q[\mbox{\boldmath$\Xi$}_0] \label{V.11}
\end{equation}
satisfies
\begin{equation}
\frac{b_0}{4a_0^2}+\frac{c_0}{a_0} >0. \label{V.12}
\end{equation}
Then 
\begin{subequations}
\begin{align}
&\sigma_0=\frac{b_0}{2a_0}+\sqrt{\frac{b_0^2}{4a_0^2}+\frac{c_0}{a_0}} \label{V.13}\\
\mbox{or}& \nonumber \\
&\sigma_0=\frac{b_0}{2a_0}-\sqrt{\frac{b_0^2}{4a_0^2}+\frac{c_0}{a_0}} 
\end{align}
\end{subequations}
In order fix the idea, suppose that \eqref{V.13} is the case. Then we can consider $\sigma$ as a functional 
of $\mbox{\boldmath$\Xi$} \in \mathfrak{H}^{\mathrm{div}}$ near $\mbox{\boldmath$\Xi$}_0$, say, $\|\mbox{\boldmath$\Xi$}-\mbox{\boldmath$\Xi$}_0\|_{\mathfrak{H}^{\mathrm{div}}}
\leq \delta ( \ll 1)$, defined by
\begin{equation}
\sigma(\mbox{\boldmath$\Xi$})=\frac{b(\mbox{\boldmath$\Xi$})}{2a(\mbox{\boldmath$\Xi$})}+
\sqrt{\frac{b(\mbox{\boldmath$\Xi$})^2}{4a(\mbox{\boldmath$\Xi$})^2}+\frac{c(\mbox{\boldmath$\Xi$})}{a(\mbox{\boldmath$\Xi$})}} \label{V.15}
\end{equation}
with
\begin{equation}
a(\mbox{\boldmath$\Xi$})=\|\mbox{\boldmath$\Xi$}\|_{\mathfrak{H}}^2,\quad
b(\mbox{\boldmath$\Xi$})=2\Omega (\mbox{\boldmath$J$}_{\star}\mbox{\boldmath$\Xi$}|\mbox{\boldmath$\Xi$})_{\mathfrak{H}},
\quad
c(\mbox{\boldmath$\Xi$})=Q[\mbox{\boldmath$\Xi$}]. \label{V.16}
\end{equation}
Here we take $\delta$ so small that 
\begin{equation}
\frac{b(\mbox{\boldmath$\Xi$})^2}{4a(\mbox{\boldmath$\Xi$})^2}+\frac{c(\mbox{\boldmath$\Xi$})}{a(\mbox{\boldmath$\Xi$})} >0
\quad\mbox{for}\quad \|\mbox{\boldmath$\Xi$}-\mbox{\boldmath$\Xi$}_0\|_{\mathfrak{H}^{\mathrm{div}}}\leq\delta.
\end{equation}\\

The variation $\mbox{\boldmath$\delta$}\sigma=\mbox{\boldmath$\delta$}\sigma(\mbox{\boldmath$\Xi$})$ of $\sigma$ at $\mbox{\boldmath$\Xi$},
\|\mbox{\boldmath$\Xi$}-\mbox{\boldmath$\Xi$}_0\|_{\mathfrak{H}^{\mathrm{div}}}\leq\delta$, is the linear functional on $\mathfrak{H}^{\mathrm{div}}$ defined by
\begin{equation}
\langle \mbox{\boldmath$\delta$}\sigma(\mbox{\boldmath$\Xi$})|\mbox{\boldmath$h$} \rangle =
\lim_{\tau \rightarrow 0}\frac{1}{\tau}
(\sigma(\mbox{\boldmath$\Xi$}+\tau \mbox{\boldmath$h$})-\sigma(\mbox{\boldmath$\Xi$})).
\end{equation}
It follows from \eqref{V.15} the equation
\begin{equation}
-\sigma^2\|\mbox{\boldmath$\Xi$}\|_{\mathfrak{H}}^2+
\sigma 2\Omega (\mbox{\boldmath$J$}_{\star}\mbox{\boldmath$\Xi$}|\mbox{\boldmath$\Xi$})_{\mathfrak{H}}+
Q[\mbox{\boldmath$\Xi$}]=0. \label{V.18}
\end{equation}
holds for $\sigma=\sigma(\mbox{\boldmath$\Xi$}), \mbox{\boldmath$\Xi$} \in \mathfrak{H}^{\mathrm{div}}, \|\mbox{\boldmath$\Xi$}-\mbox{\boldmath$\Xi$}_0\|_{\mathfrak{H}^{\mathrm{div}}}\leq\delta$.
Therefore we have
\begin{align*}
&(-2\sigma \|\mbox{\boldmath$\Xi$}\|_{\mathfrak{H}}^2+2\Omega (\mbox{\boldmath$J$}_{\star}\mbox{\boldmath$\Xi$}|\mbox{\boldmath$\Xi$})_{\mathfrak{H}})\mbox{\boldmath$\delta$}\sigma + \\
&-\sigma^2\mbox{\boldmath$\delta$}\|\mbox{\boldmath$\Xi$}\|_{\mathfrak{H}}^2
+\sigma \mbox{\boldmath$\delta$}(2\Omega J_{\star}\mbox{\boldmath$\Xi$}|\mbox{\boldmath$\Xi$})_{\mathfrak{H}}
+\mbox{\boldmath$\delta$} Q[\mbox{\boldmath$\Xi$}] =0,
\end{align*}
or, precisely writing,
\begin{align*}
&(-2\sigma \|\mbox{\boldmath$\Xi$}\|_{\mathfrak{H}}^2+2\Omega (\mbox{\boldmath$J$}_{\star}\mbox{\boldmath$\Xi$}|\mbox{\boldmath$\Xi$})_{\mathfrak{H}})\langle \mbox{\boldmath$\delta$}\sigma(\mbox{\boldmath$\Xi$})|\mbox{\boldmath$\delta$}\mbox{\boldmath$\Xi$}\rangle + \\
&
+2\mathfrak{Re}\Big[
-\sigma^2
(\mbox{\boldmath$\Xi$} |\mbox{\boldmath$\delta$}\mbox{\boldmath$\Xi$})_{\mathfrak{H}}
+\sigma 2\Omega (J_{\star}\mbox{\boldmath$\Xi$}|\mbox{\boldmath$\delta$} \mbox{\boldmath$\Xi$})_{\mathfrak{H}}
+ Q(\mbox{\boldmath$\Xi$}, \mbox{\boldmath$\delta$}\mbox{\boldmath$\Xi$})
\Big] =0
\end{align*}
for $\forall \mbox{\boldmath$\delta$}\mbox{\boldmath$\Xi$} \in \mathfrak{H}^{\mathrm{div}}$.  
Here we note that 
\begin{align*}
-2\sigma \|\mbox{\boldmath$\Xi$}\|_{\mathfrak{H}}^2+2\Omega (\mbox{\boldmath$J$}_{\star}\mbox{\boldmath$\Xi$}|\mbox{\boldmath$\Xi$})_{\mathfrak{H}}&=
-2\sigma a(\mbox{\boldmath$\Xi$})+b(\mbox{\boldmath$\Xi$}) = \\
&=-\sqrt{
\frac{b(\mbox{\boldmath$\Xi$})^2}{4a(\mbox{\boldmath$\Xi$})^2}+\frac{c(\mbox{\boldmath$\Xi$})}{a(\mbox{\boldmath$\Xi$})} 
}\not=0 
\end{align*}
for $\|\mbox{\boldmath$\Xi$}-\mbox{\boldmath$\Xi$}_0\|_{\mathfrak{H}^{\mathrm{div}}}\leq \delta$.
Therefore $\mbox{\boldmath$\delta$} \sigma(\mbox{\boldmath$\Xi$})=0$ if and only if
$$-\sigma^2
(\mbox{\boldmath$\Xi$} |\mbox{\boldmath$h$})_{\mathfrak{H}}
+\sigma 2\Omega (J_{\star}\mbox{\boldmath$\Xi$}|\mbox{\boldmath$h$})_{\mathfrak{H}}
+ Q(\mbox{\boldmath$\Xi$}, \mbox{\boldmath$h$})
=0$$
for $\forall \mbox{\boldmath$h$} \in \mathfrak{H}^{\mathrm{div}}$.
Thus we can claim the following `variational principle':

\begin{Theorem}
Let $\mbox{\boldmath$\Xi$}_0 \in \mathfrak{H}^{\mathrm{div}}$ satisfy 
$$\|\mbox{\boldmath$\Xi$}_0\|_{\mathfrak{H}}^2 >0, \quad
\frac{(2\Omega
 (\mbox{\boldmath$J$}_{\star}\mbox{\boldmath$\Xi$}_0|\mbox{\boldmath$\Xi$}_0)_{\mathfrak{H}})^2}{4\|\mbox{\boldmath$\Xi$}_0\|_{\mathfrak{H}}^4}+
\frac{Q[\mbox{\boldmath$\Xi$}_0]}{\|\mbox{\boldmath$\Xi$}_0\|_{\mathfrak{H}}^2} >0.
$$
The variation $\mbox{\boldmath$\delta$} \sigma$ of $\sigma$ (specified by \eqref{V.15}) vanishes at $\mbox{\boldmath$\Xi$}_0$ if and only if 
$\mbox{\boldmath$\Xi$}_0 \in \mathsf{D}(\mbox{\boldmath$L$})$ and enjoys the equation \eqref{EigEq}:
$$
-\sigma_0^2\mbox{\boldmath$\Xi$}_0+\sigma_0 2\Omega \mbox{\boldmath$J$}_{\star}\mbox{\boldmath$\Xi$}_0+\mbox{\boldmath$L$}\mbox{\boldmath$\Xi$}_0=\mathbf{0}.
$$
Then $e^{\mathrm{i}\sigma_0t}\mbox{\boldmath$\Xi$}_0(\mbox{\boldmath$x$})$ is a solution of the equation \eqref{V.3}. Here $\sigma_0=\sigma(\mbox{\boldmath$\Xi$}_0)$.
\end{Theorem}

This principle tells us that, if we want to find an eigenfrequency, we may try to find a stationary point of the functional
$$\sigma(\mbox{\boldmath$\Xi$})=\frac{b(\mbox{\boldmath$\Xi$})}{2}\pm\sqrt{\frac{b(\mbox{\boldmath$\Xi$})^2}{4}+c(\mbox{\boldmath$\Xi$})}
$$ under the constraint $\|\mbox{\boldmath$\Xi$}\|_{\mathfrak{H}}=1$. But it seems that this principle is far from
the solution of the problem to establish the existence and completeness of the system of eigenvectors.

 For example, as a `practical use of the variational principle', D. Lynden-Bell and J. P. Ostriker, \cite[Section 2.6]{LyndenBO}, proposed the following scheme:

\begin{quote}

Take sufficiently many functions $\mbox{\boldmath$\Xi$}_{(i)}, i=1,\cdots, N$ as those who consist a base, and consider the trial function
$$\mbox{\boldmath$\Xi$}=\sum_{i=1}^Na^i\mbox{\boldmath$\Xi$}_{(i)},\quad
\mbox{\boldmath$a$} =
\begin{bmatrix}
a^1 \\
\cdot \\
\cdot \\
a^N
\end{bmatrix}.
$$
Put
\begin{align*}
& \mbox{\boldmath$A$}= (A_{ij})_{i,j},\quad A_{ij}=(\mbox{\boldmath$\Xi$}_{(i)}|\mbox{\boldmath$\Xi$}_{(j)})_{\mathfrak{H}}, \\
& \mbox{\boldmath$B$}=(B_{ij})_{i,j} ,\quad B_{ij}=2\Omega(\mbox{\boldmath$J$}_{\star}\mbox{\boldmath$\Xi$}_{(i)}|\mbox{\boldmath$\Xi$}_{(j)})_{\mathfrak{H}}, \\
& \mbox{\boldmath$C$}=(C_{ij})_{i,j} , \quad C_{ij}=Q(\mbox{\boldmath$\Xi$}_{(i)}, \mbox{\boldmath$\Xi$}_{(j)}).
\end{align*}
Then $$\sigma(\mbox{\boldmath$\Xi$})=
\Big((-\sigma^2\mbox{\boldmath$A$}+\sigma \mbox{\boldmath$B$}+\mbox{\boldmath$C$})\mbox{\boldmath$a$} \Big| \mbox{\boldmath$a$}\Big).
$$
Thus, D. Lynden-Bell and J. P. Ostriker say, the variational principle reads 
$$\mbox{\boldmath$\delta$} \Big((-\sigma^2\mbox{\boldmath$A$}+\sigma \mbox{\boldmath$B$}+\mbox{\boldmath$C$})\mbox{\boldmath$a$}\Big|\mbox{\boldmath$a$}\Big)=0
$$
and, varying $\mbox{\boldmath$a$}$, we obtain the `secular determinant'
\begin{equation}
\mathrm{det}
(-\sigma^2\mbox{\boldmath$A$}+\sigma \mbox{\boldmath$B$}+\mbox{\boldmath$C$} ) =0\label{V.det}
\end{equation}
for the determination of the variationally best eigenfrequencies $\sigma$.
\end{quote}
But it is doubtful that this scheme is so practtical, since we have no confidence that the equation \eqref{V.det}, which is an algebraic equation for the unknown $\sigma$ of degree is $2N$, can be numerically solved to determine an approximating eigenfrequency $\sigma_N$ so that they converge to a true eigenfrequency as $N \rightarrow \infty$. 

On the other hand, since
$$\sigma(\mbox{\boldmath$\Xi$})=\frac{b(\mbox{\boldmath$\Xi$})}{2}+
\sqrt{\frac{b(\mbox{\boldmath$\Xi$})^2}{4}+c(\mbox{\boldmath$\Xi$})} \geq 0
$$
for $\forall \mbox{\boldmath$\Xi$} \in \mathfrak{H}^{\mathrm{div}}$, we can consider 
$$\sigma_*:=\inf\{ \sigma(\mbox{\boldmath$\Xi$}) | \mbox{\boldmath$\Xi$} \in \mathfrak{H}^{\mathrm{div}}, \|\mbox{\boldmath$\Xi$}\|_{\mathfrak{H}}=1 \}
$$ and we may expect that the minimum might give an eigenfrequency. But, when $\Omega \not=0$, we are not sure about the existence of a $\mbox{\boldmath$\Xi$}_*$  which attains the infimum, namely $ \sigma_*=\sigma(\mbox{\boldmath$\Xi$}_*)$,  in general, since the imbedding
$\mathfrak{H}^{\mathrm{div}} \hookrightarrow \mathfrak{H}$ is not compact
and we may be unable to extract convergent subsequences from a minimizing sequence, say, a sequence
$(\bm{\Xi}_n)_{n\in\mathbb{N}}$ such that 
$\bm{\Xi}_n \in \mathfrak{H}^{\mathrm{div}}, \|\bm{\Xi}_n\|_{\mathfrak{H}}=1,
\sigma(\bm{\Xi}_n) \rightarrow \sigma_*$ as $n \rightarrow \infty$.

 When $\Omega=0$, then $ b(\mbox{\boldmath$\Xi$})=0$, and the variational problem
$$\lambda_*(=\sigma_*^2)=\inf\{ c(\mbox{\boldmath$\Xi$}) (=Q[\mbox{\boldmath$\Xi$}] ) | \mbox{\boldmath$\Xi$} \in \mathfrak{H}^{\mathrm{div}}, \|\mbox{\boldmath$\Xi$}\|_{\mathfrak{H}}=1 \}
$$
actually admits the trivial solution $\lambda_*=0$ with any eigenvector $\in \mathrm{Ker}(\mbox{\boldmath$L$})$.
But the Mini-Max principle does not work, for, since
$\mathrm{dim.Ker}(\mbox{\boldmath$L$})=\infty$, we cannot go ahead across the $0$ eigenvalue eternally, although there actually  remain infinitely many positive eigenvalues.
Thus also in this case the variational principle seems to be not so useful.

However, when $\Omega=0$, the equation \eqref{EigEq} reduces to
\begin{equation}
-\lambda\bm{\Xi}+\bm{L}\bm{\Xi}=\mathbf{0}, \label{501}
\end{equation}
where $\lambda=\sigma^2$, and this eigenvalue problem is completely solved as follows: \\

Suppose $\Omega=0$. We note that the background stationary solution $\bar{\rho}$ is a spherically symmetric equilibrium, say,
\begin{equation}
\bar{\rho}(\bm{x})=\Big(\frac{(\gamma-1)\mathsf{G}M_0}{\mathsf{A}\gamma}\Big(\frac{1}{r}-\frac{1}{R}\Big)\Big)^{\frac{1}{\gamma-1}} \label{502}
\end{equation}
and $\mathfrak{R}=\{\bar{\rho} >0\}=\{ R_0 < r < R\}$ is an annulus. Taking the divergence of $\bar{\rho}$ times \eqref{501}, the problem reduces to
\begin{equation}
-\lambda g-\mathrm{div}\bar{\rho}\Big(\mathrm{grad}
\overline{\frac{d\varUpsilon}{d\rho}} g\Big)=0, \label{503}
\end{equation}
where 
\begin{equation}
g=\mathrm{div}(\bar{\rho}\bm{\Xi}). \label{504}
\end{equation}
Note that we can treat the operator $\mathcal{N}: g \mapsto
\mathrm{div}\bar{\rho}\Big(\mathrm{grad}
\overline{\frac{d\varUpsilon}{d\rho}} g\Big)$ as  that with the similar property to the Laplacian operator $\triangle=\mathrm{div.grad}$, taking into account the singular behaviors of $\bar{\rho}$ and
$\displaystyle \overline{\frac{d\varUpsilon}{d\rho}}$ at the physical vacuum boundary, say,
\begin{align*}
&\bar{\rho} \sim \Big(\frac{(\gamma-1)\mathsf{G}M_0}{\mathsf{A}\gamma}\Big)^{\frac{1}{\gamma-1}}
(R-r)^{\frac{1}{\gamma-1}}, \\
&\overline{\frac{d\varUpsilon}{d\rho}} \sim 
\mathsf{A}\gamma \Big(\frac{(\gamma-1)\mathsf{G}M_0}{\mathsf{A}\gamma}\Big)^{\frac{\gamma-2}{\gamma-1}}
(R-r)^{-\frac{2-\gamma}{\gamma-1}} 
\end{align*} 
near $\Sigma_1=\{ r= R\}$. Therefore we can claim:\\

{\it The operator $\mathcal{N}$ can be considered as a self-adjoint operator in the Hilbert space $\mathfrak{G}$
 and its spectrum $\sigma(\mathcal{N})$ is of the Sturm-Liouville type, that is, 
$\sigma(\mathcal{N})=\{ \lambda_n^{\mathcal{N}} | n \in \mathbb{N}\}$, where $\lambda_n^{\mathcal{N}}$ is an eigenvalue with finite multiplicity,
$0<\lambda_0^{\mathcal{N}}<\cdots <
\lambda_n^{\mathcal{N}} < \lambda_{n+1}^{\mathcal{N}}$,  and $\lambda_n^{\mathcal{N}} \rightarrow +\infty $ as $ n \rightarrow \infty$.
Here $$\mathfrak{G}=
 L^2\Big(\mathfrak{R}; \overline{\frac{d\varUpsilon}{d\rho}}d\bm{x}\Big)
\cap\Big\{ g \  \Big|\  \int_{\mathfrak{R}}gd\bm{x}=0\Big\}.$$
}\\

Note that the imbedding of $ L^2(\mathfrak{R}; \overline{\frac{d\varUpsilon}{d\rho}}d\bm{x})$ into $L^2(\mathfrak{R}; d\bm{x})
(\hookrightarrow
L^1(\mathfrak{R};d\bm{x}))$ is continuous. 
For a proof see the proof of \cite[Theorem 2]{JJTM2020}. Hence the argument of \cite{JJTM2020} leads us to the following conclusion:\\

{\it When $\Omega=0$, $\bm{L}$ can be considered as a self-adjoint operator in the Hilbert space $\mathfrak{F}$ and its spectrum $\sigma(\bm{L})$ coincides with
$\sigma(\mathcal{N}) \cup \{0\}$, while $\mathrm{dim.Ker}(\bm{L})=\infty$ and  $\lambda=\lambda_n^{\mathcal{N}}\not=0$ is an eigenvalue with finite multilicity.
Here
$$\mathfrak{F}=
\{ \bm{\Xi}\in\mathfrak{H}\  |\  \mathrm{div}(\bar{\rho}\bm{\Xi} )\in \mathfrak{G}\}
=\Big\{ \bm{\Xi}\in \mathfrak{H}^{\mathrm{div}}\ \Big|\  
\int_{\mathfrak{R}} \mathrm{div}(\bar{\rho}\bm{\Xi})d\bm{x}=0\Big\}. $$
}\\

Note that $\mathfrak{F}$ is a closed subspace of $\mathfrak{H}^{\mathrm{div}}$, since
$ \mathfrak{G} \ni g \mapsto \int_{\mathfrak{R}}gd\bm{x}$ is constinuous. 
In fact, we have $\displaystyle \int_{\mathfrak{R}}\overline{\frac {d\rho}{d\Upsilon}}d\bm{x}
<+\infty$ thanks to $\gamma >1$ so that
$$\Big| \int g \Big| \leq
\Big[\int\overline{\frac{d\Upsilon}{d\rho}}|g|^2\Big]^{1/2}
\Big[\int_{\mathfrak{R}}\overline{\frac{ d\rho}{d\Upsilon}} \Big]^{1/2}
\lesssim \|g\|_{L^2(\frac{d\Upsilon}{d\rho}d\bm{x})}.$$  
We have $\mathfrak{H}_0^{\mathrm{div}} \subset \mathfrak{F}$, since
$$\int \mathrm{div}(\bar{\rho}\bm{\varphi})d\bm{x}=0 \quad\mbox{for}\quad 
\forall \bm{\phi} \in C_0^{\infty}(\mathfrak{R}).$$

Therefore by the well-known theorem (
\cite[p.905, X.3.4.Theorem]{DunfordS} complemented by 
\cite[p. 177, Chapter III, Theorem 6.15]{Kato}), 
we can say that the eigenvectors of a CONS of $\mathrm{Ker}(\bm{L})$
and all
$$\psi_n=
\frac{1}{\lambda_n^{\mathcal{N}}}\bar{\rho}\mathrm{grad}\Big(-
\overline{\frac{d\varUpsilon}{d\rho}}\varphi_n^{\mathcal{N}}\Big)
\Big\|
\frac{1}{\lambda_n^{\mathcal{N}}}\bar{\rho}\mathrm{grad}\Big(-
\overline{\frac{d\varUpsilon}{d\rho}}\varphi_n^{\mathcal{N}}\Big)
\Big\|^{-1},
$$
$\varphi_n^{\mathcal{N}}$ being an eigenvalue of $\mathcal{N}$ associated with the eigenvalue $\lambda_n^{\mathcal{N}}\not=0$ form a complete ortogonal system of the Hilbert space
 $\mathfrak{F}$. \\
 
 In this sense, when $\Omega=0$, the eigenfrequency problem is completely solved.

\begin{Remark}
Let $\Omega=0$.
If we consider the operator $\bm{L}$ in the space $\mathfrak{H}$, we can claim that $\{0\}\cap \sigma(\mathcal{N})=\{ 0, \lambda_0^{\mathcal{N}}, \lambda_1^{\mathcal{N}},\cdots \} \subset
\sigma(\bm{L})$, but we do not know whether they coincide or not, say,
we do not know whether there are real continuous spectrum between the eigenvalues or not. 
\end{Remark}

However, when $\Omega\not=0$, the above discussion seems not to work. In fact the term
$2\Omega(\nabla|\bar{\rho}J_{\star}\bm{\Xi})$
may cause trouble, since it cannot be reduced to a quantity determined by $g=(\nabla|\bar{\rho}\bm{\Xi})$.


\vspace{15mm}

{\bf\Large Acknowledgment}\\

The author expresses his sincere thanks to Professor Akitaka Matsumura (Osaka University) for kindly having helpful discussions on this research and especially enlghting the author on the boundary trace of the divergence spaces $H(\mathrm{div}), H_0(\mathrm{div})$. This work is supported by JSPS KAKENHI Grant Number JP21K03311.\\

{\bf\Large  The data availability statement}\\

No new data were created or analyzed in this study.


\end{document}